\documentclass[final,3p]{elsarticle}

\usepackage[colorlinks=true,breaklinks=true,pdftex]{hyperref}

\journal{CAGD 2021}

\biboptions{authoryear}

\usepackage{graphics,graphicx,amssymb,amstext,amsmath,amsfonts}
\usepackage{wrapfig,multicol}
\DeclareMathOperator{\W}{\mathsf{W}}
\DeclareMathOperator{\fv}{\mathsf{fv}}
\DeclareMathOperator{\B}{\mathsf{B}}
\DeclareMathOperator{\R}{\mathsf{R}}
\DeclareMathOperator{\M}{\mathsf{M}}
\DeclareMathOperator{\A}{\mathsf{A}}
\DeclareMathOperator{\Le}{\mathsf{L}}
\DeclareMathOperator{\I}{\mathbf{i}}

\newtheorem{myLemma}{Lemma}[section]
\newtheorem{myDef}{Definition}[section]

\begin{document}

\begin{frontmatter}

\title{A Simple and Complete Discrete Exterior Calculus on General Polygonal Meshes} 

\author{Lenka Pt\'a\v ckov\'a\footnote{Corresponding author, email: \url{lenaptackova@gmail.com}}, Luiz Velho}


%
%
\begin{abstract}
 Discrete exterior calculus (DEC) offers a coordinate--free discretization of exterior calculus especially suited for computations on curved spaces. In this work, we present an extended version of DEC on surface meshes formed by general polygons that bypasses the need for combinatorial subdivision and does not involve any dual mesh. At its core, our approach introduces a new polygonal wedge product that is compatible with the discrete exterior derivative in the sense that it satisfies the Leibniz product rule. Based on the discrete wedge product, we then derive a novel primal--to--primal Hodge star operator. Combining these three `basic operators' we then define new discrete versions of the contraction operator and Lie derivative, codifferential and Laplace operator. We discuss the numerical convergence of each one of these proposed operators and compare them to existing DEC methods. Finally, we show simple applications of our operators on Helmholtz--Hodge decomposition, Laplacian surface fairing, and Lie advection of functions and vector fields on meshes formed by general polygons.
\end{abstract}

\begin{keyword}
Discrete exterior calculus \sep Geometry processing with polygonal meshes \sep Polygonal wedge product \sep Polygonal Hodge star operator \sep Discrete Laplace operator \sep Discrete Lie advection
\end{keyword}

\end{frontmatter}


\begin{figure}[h]
\centering{
\includegraphics[width=0.24\linewidth]{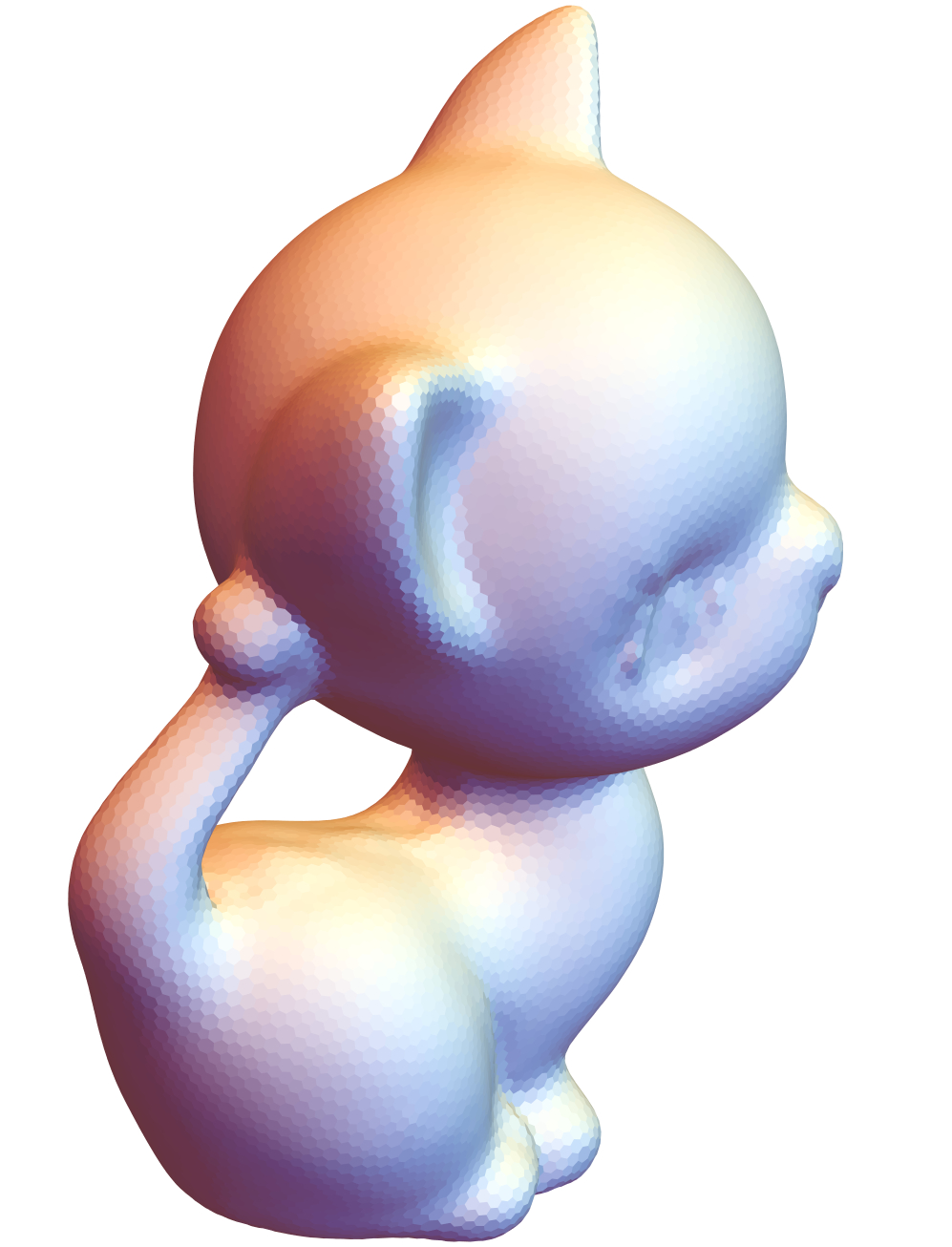}
\includegraphics[width=0.24\linewidth]{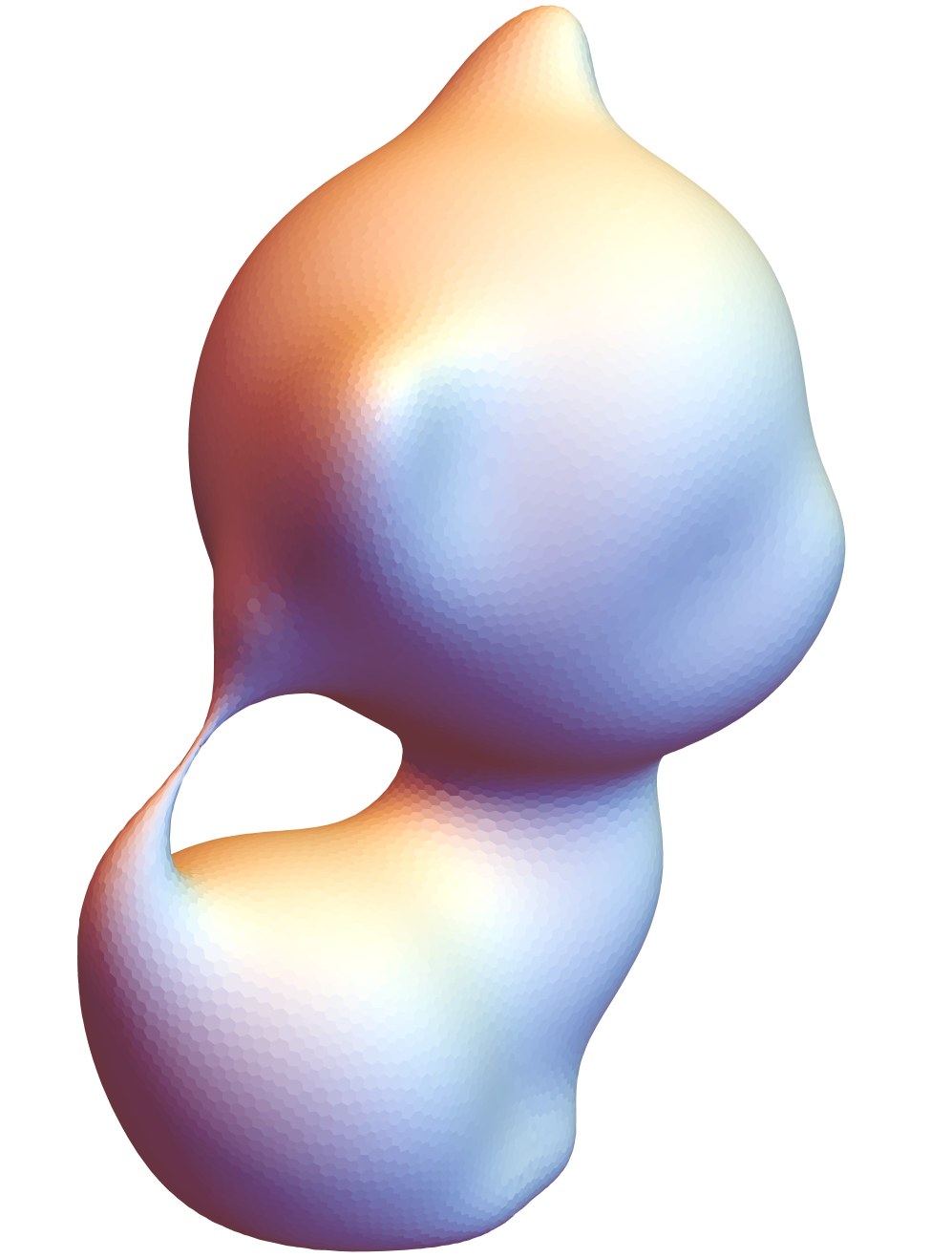}
\includegraphics[width=0.24\linewidth]{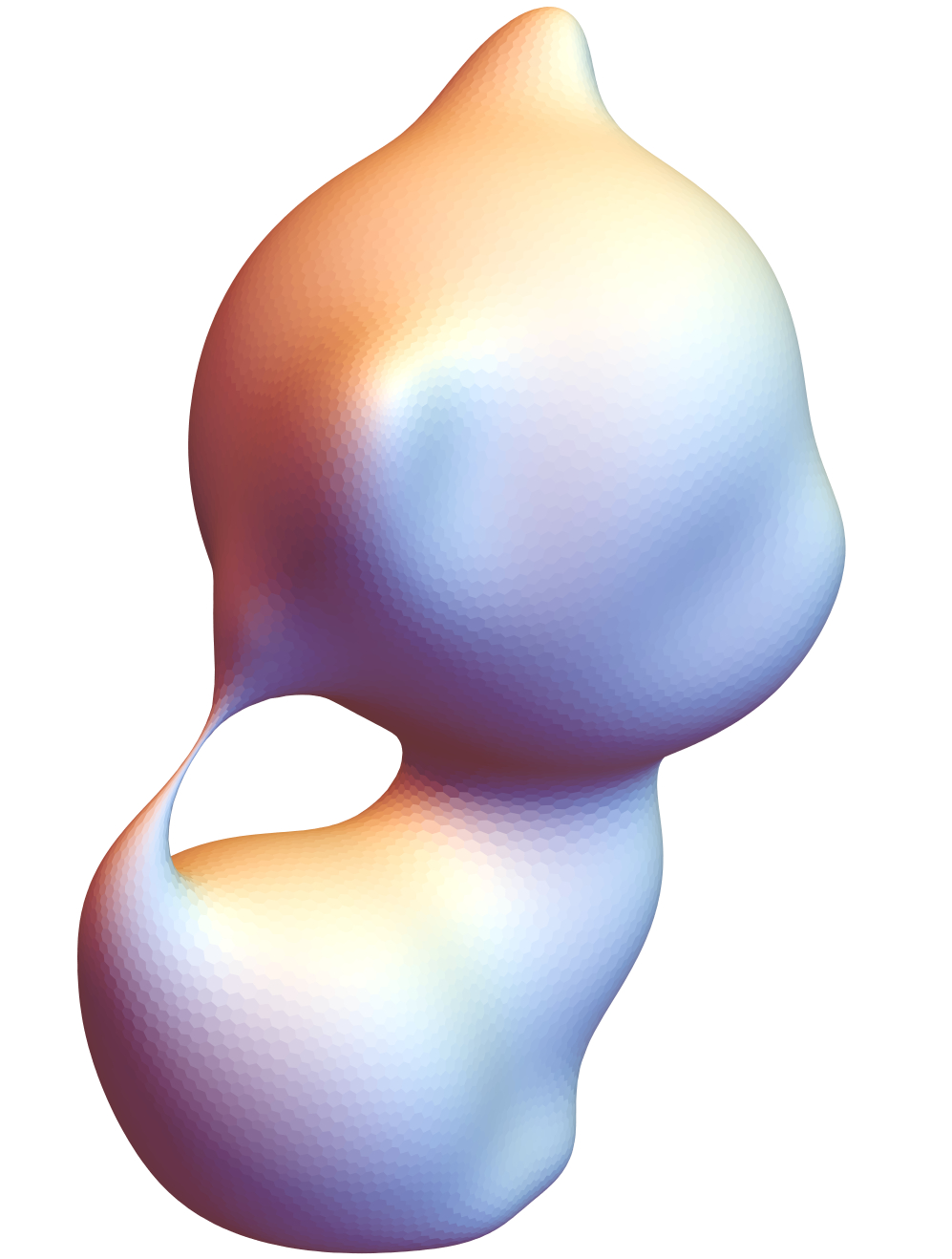}
\includegraphics[width=0.24\linewidth]{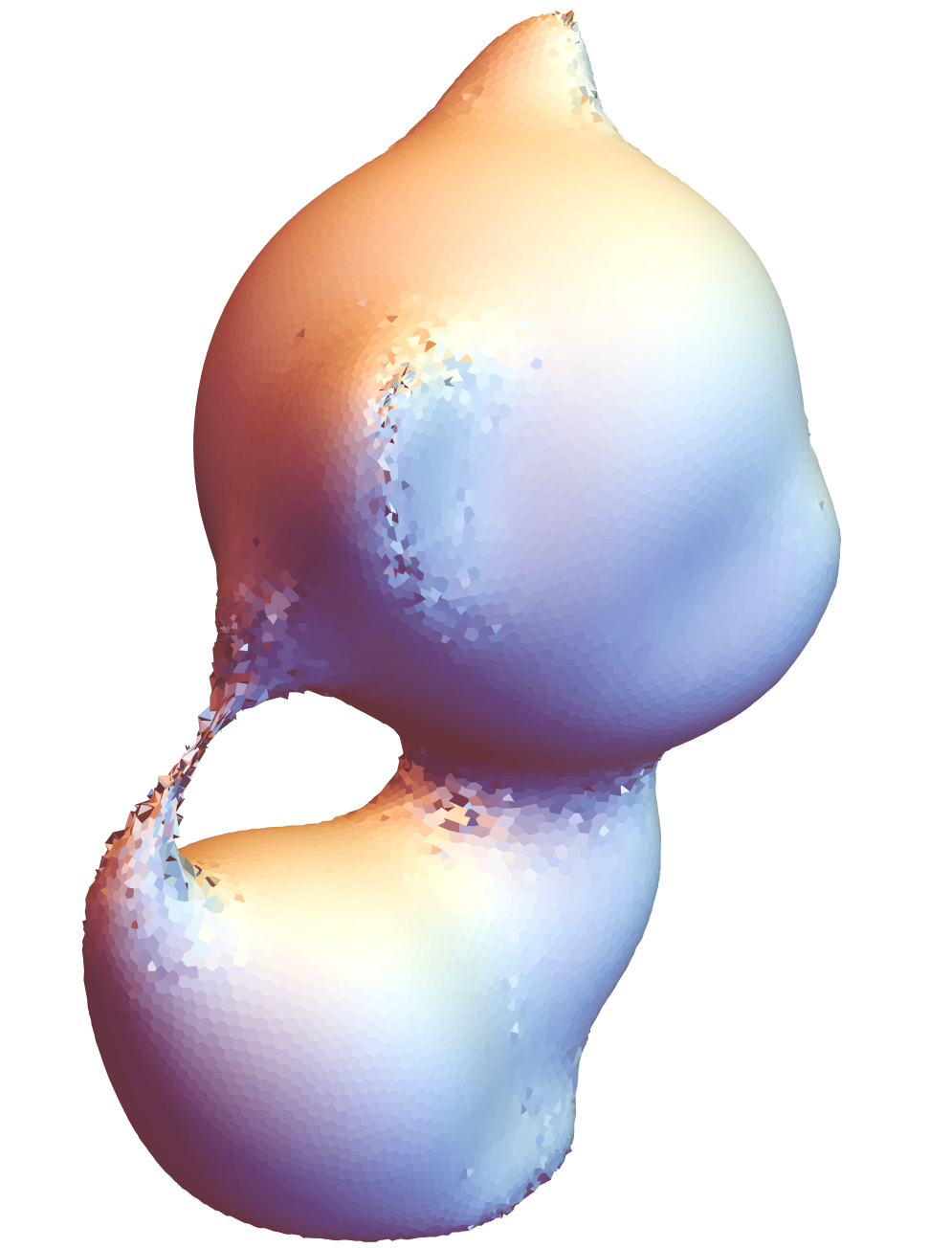}}
\caption{Comparison of implicit mean curvature flows on a general polygonal mesh (29k vertices) after 10 iterations with time step $t = 10^{-4}$. On the far left is the original mesh.
Our method (center left) and the algorithm of \citet{Alexa2011} with their combinatorially enhanced Laplacian (center right) produce visually well--smoothed meshes. However, their method with purely geometric Laplacian (far right) exhibits some undesirable artifacts on the ears, neck, and tail of the kitten.}
\label{fig:LaplaceFlowKitty}
\end{figure}

\section{Introduction}
The discretization of differential operators on surfaces is fundamental for geometry processing tasks, ranging from remeshing to vector fields manipulation. Discrete exterior calculus (DEC) is arguably one of the prevalent numerical frameworks to derive such discrete differential operators. However, the vast majority of work on DEC is restricted to simplicial meshes, and far less attention has been given to meshes formed by arbitrary polygons, possibly non--planar and non--convex.

In this work, we propose a new discretization for several operators commonly associated to DEC that operate directly on polygons without involving any subdivision.
Our approach offers three main practical benefits. First, by working directly with polygonal meshes, we overcome the ambiguities of subdividing a discrete surface into a triangle mesh. Second, our construction operates solely on primal elements, thus removing any dependency on dual meshes. Finally, our method includes the discretization of new differential operators such as the contraction operator and Lie derivatives.

We concisely expose our framework, describe each of our operators and compare them to existing DEC methods.
We examine the accuracy of our numerical scheme by a series of convergence tests on flat and curved surface meshes. We also demonstrate the applicability of our method for Helmholtz--Hodge decomposition of vector fields, surface fairing, and Lie advection of vector fields and functions.



\section{Related Work and Preliminaries}\label{sec:RelatedWork}
There is a vast literature on DEC on triangle meshes, e.g., \citep{Hirani,DEC2005,DDFormsForCM,SiggraphCourse} -- all these publications have in common that they deal with purely simplicial meshes and use a dual mesh to define operators, we will refer to their approach as to the \textbf{classical DEC}.

As announced, unlike the classical DEC, our method works with general polygonal meshes and does not involve any dual meshes. However, our operators differ also in other aspects, e.g., support (see Figures \ref{fig:codifferential1}--\ref{fig:triLaplacian0}). Next we briefly introduce several basic DEC notions and point out the key differences between our approach and existing schemes, principally the classical DEC.

\subsection{Discrete differential forms and the exterior derivative}
We strictly stick to the convention, common to previous DEC literature, that a discrete $q$--form is located on $q$--dimensional cells of the given mesh.

Discrete differential forms are usually denoted by small Greek letters and sometimes we add a number superscript to emphasize the degree of the form, i.e., a $q$--form $\alpha$ can be denoted as $\alpha^q$. 

A polygonal mesh $S$ is made of a set of vertices (0--dimensional cells), oriented edges (1--dimensional cells), and oriented faces (2--dimensional cells). A real discrete differential $q$--form $\alpha^q$ on $S$ is a $q$--cochain, i.e., a real number assigned to each $q$--dimensional cell $c^q$ of $S$. E.g., if $(e_0,e_1,\dots,e_n)$ is the vector of all edges of $S$, then an 1--form $\alpha^1$ is a vector of real values
\[
 \alpha^1=(\alpha(e_0),\dots,\alpha(e_n)).
\]
The \textbf{discrete exterior derivative} $d$ is the coboundary operator and it holds:
\[
  (d\alpha)(c^{q+1}) = \alpha(\partial c^{q+1}) = \sum_{c^q\in S} [c^{q+1}:c^q]\alpha(c^q),
\]
where $\partial$ is the boundary operator and $[c^{q+1}:c^q]$ denotes the incidence relation between cells $c^{q+1}$ and $c^{q}$, as depicted in the example bellow.

\begin{minipage}[c]{0.15\linewidth}
\includegraphics[width=\linewidth]{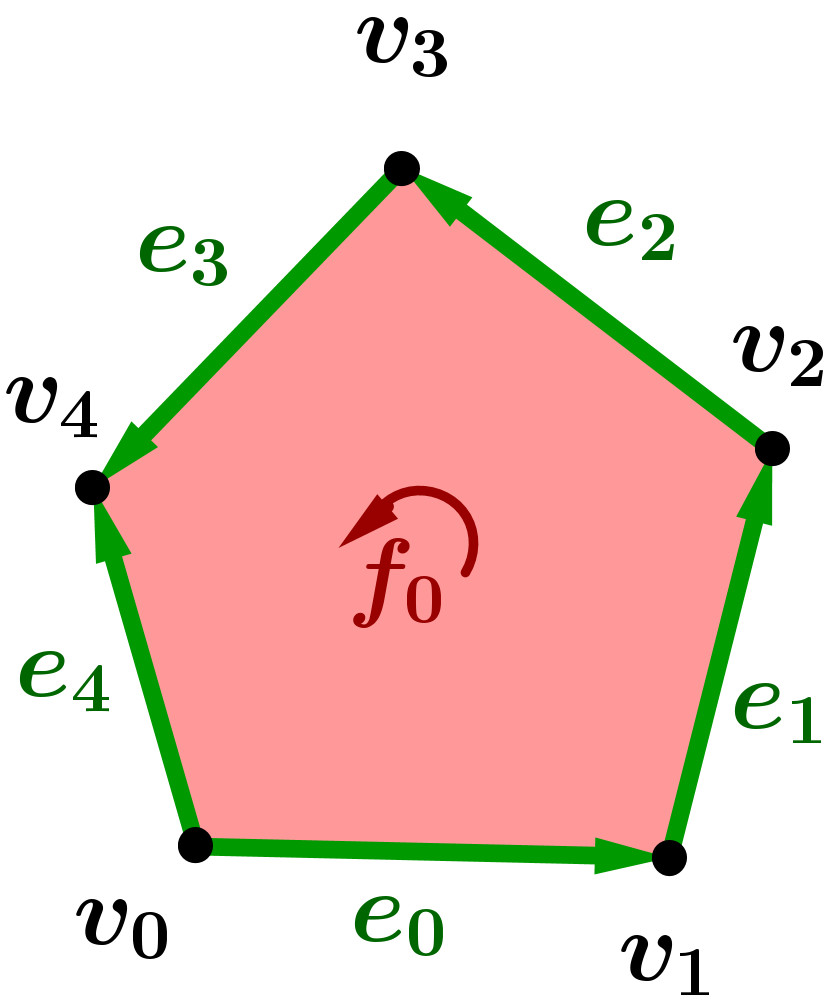}
\end{minipage}
\begin{minipage}[l]{0.8\linewidth}
\small
\begin{align*}
[f_0:e_0] &=1,[f_0:e_1]=1,\dots, [f_0:e_4]=-1 \\
\partial f_0 & =\sum\limits_{e_i\in S} [f_0:e_i]e_i \;=\; e_0+e_1 +e_2+e_3-e_4 \\
\alpha^1 &= (\alpha(e_0),\dots,\alpha(e_n)) \\
d\alpha(f_0) &= \sum\limits_{e_i\in S}[f_0:e_i]\alpha(e_i) \;=\; \alpha(e_0)+\alpha(e_1) +\alpha(e_2)+\alpha(e_3) -\alpha(e_4) 
\end{align*}
\normalsize
\end{minipage}
\vspace{0.4pt}

\noindent
The boundary of a face $f_0$ is a sum of incident oriented edges, where we take in account the orientation of the boundary edges with respect to the given face. The discrete exterior derivative of a 1--form $\alpha^1$ (stored on edges) is a 2--form $d\alpha$ located on faces and it is the ``oriented sum'' of the values of $\alpha$ on boundary edges of $f_0$.

\subsection{The cup product and the wedge product}
We consider the wedge product on polygons to be the main building block of our theory. On smooth manifolds, the wedge product allows for building higher degree forms from lower degree ones. Similarly in algebraic topology of pseudomanifolds, a cup product is a product of two cochains of arbitrary degree $p$ and $q$ that returns a cochain of degree $p+q$ located on $(p+q)$--dimensional cells. Thus we consider the cup product to be the appropriate discrete version of the wedge product.

The cup product was introduced by J. W. Alexander, E. \v Cech, and H. Whitney in 1930's and it became a well--studied notion in algebraic topology, mainly in the simplicial setting. Later, the cup product was extended also to $n$--cubes~\citep{Massey,Rachel}.

In graphics, \citet{GlobalConformalParametrization} presented a wedge product of two discrete 1--forms on triangulations, which turns to be equivalent to the cup product of two 1--cochains on triangle complexes well studied in algebraic topology, e.g., \citep{Whitney}.

\begin{figure}[ht]
\centering\includegraphics[width=0.7\linewidth]{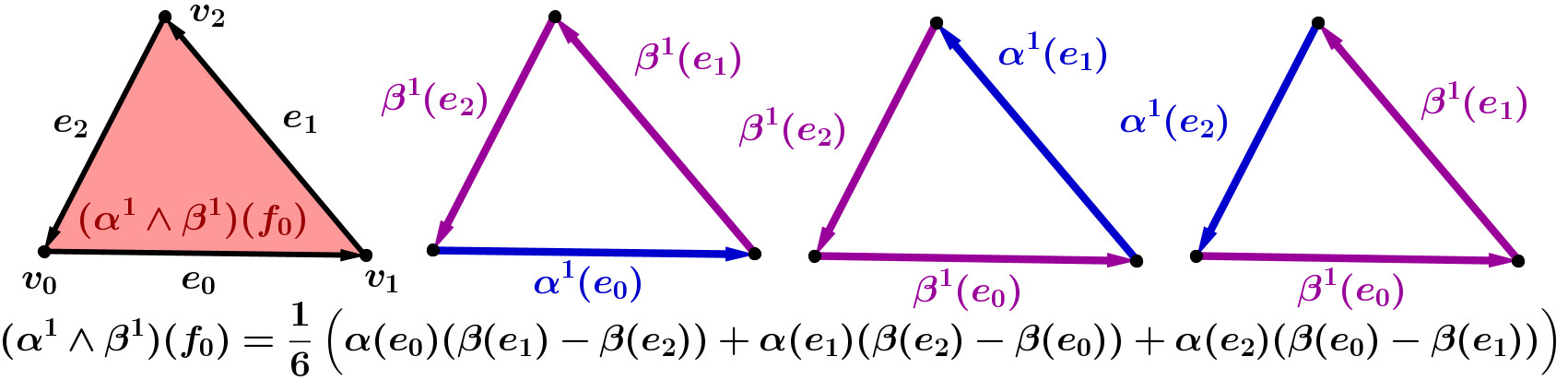}
\caption{The wedge product of two 1--forms on a triangle: the product of two 1--forms is a 2--form located on faces (far left). }\label{fig:CupTriangle11}
\end{figure}

On triangles our discrete wedge product is equivalent to the cup product of \citet{Whitney}, see also the Figure \ref{fig:CupTriangle11}. On quadrilaterals our discrete wedge product is equivalent to the cup product of \citet{Massey}.

In common to previous approaches, the discrete wedge product is metric--independent and satisfies core properties such as the Leibniz product rule, skew--commutativity, and associativity on closed forms.

\subsection{The Hodge star operator}
The most common discretization of the Hodge star operator on triangle meshes is the so called diagonal approximation, see, e.g., \citep{DDFormsForCM}, which is computed based on the ratios between the volumes of primal simplices and their dual cells.
In contrast, we propose a Hodge star operator that does not use a dual mesh. Since our dual forms are again located on primal elements, we can compute the wedge product of primal and dual forms and hence define further operators.    
However, this primal--primal definition brings some drawbacks as well, we discuss them in Section \ref{subsec:HodgeStar}.

\subsection{The Hodge inner product}
On smooth manifolds the Hodge star operator together with the wedge product define the Hodge inner product, our definition is derived in the same fashion. 

On the contrary, in classical DEC, e.g, \citep{DDFormsForCM}, the Hodge star is actually derived from a previously given inner product. \citet{Alexa2011} in Lemma 3 also present a discrete version of inner product matrices, that become the building blocks of their theory.

\subsection{The codifferential}
On a Riemannian $n$--manifold, the Hodge star operator is employed to define the codifferential operator $\delta(\alpha^k)=(-1)^{n(k-1)+1} \star d\star\alpha.$
It is a linear operator that maps $k$--forms to $(k-1)$-- forms. On 1--forms it is also called the divergence operator.

\begin{figure}[htb]
\centering\includegraphics[width=0.4\linewidth]{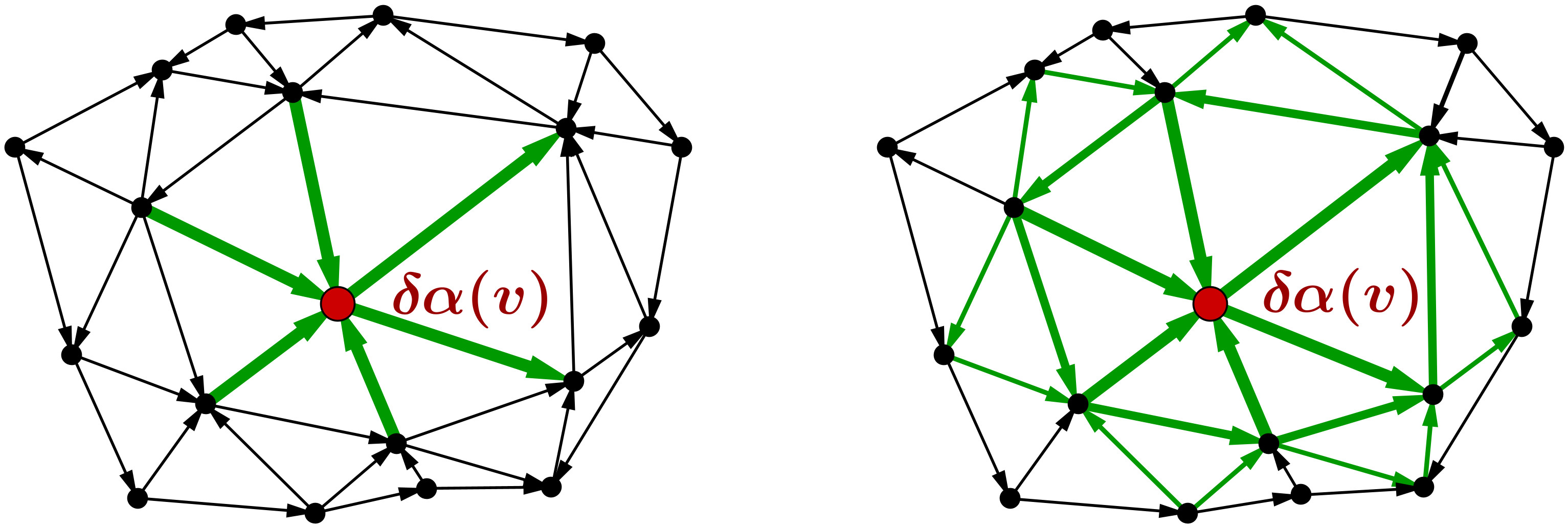}
\caption{Comparison of the support of the codifferential of 1--forms between the classical DEC (L) and our method (R). The codifferential of a 1--form $\alpha$ is a 0--form located on vertices. The value of $\delta\alpha$ on the red vertex $v$ is a linear combination of values of $\alpha$ on edges colored green. The edge thickness reflects the weight of the corresponding edge values $\alpha$ on $\delta\alpha(v)$.
}\label{fig:codifferential1}
\end{figure}

\begin{figure}[htb]
\centering\includegraphics[width=0.42\linewidth]{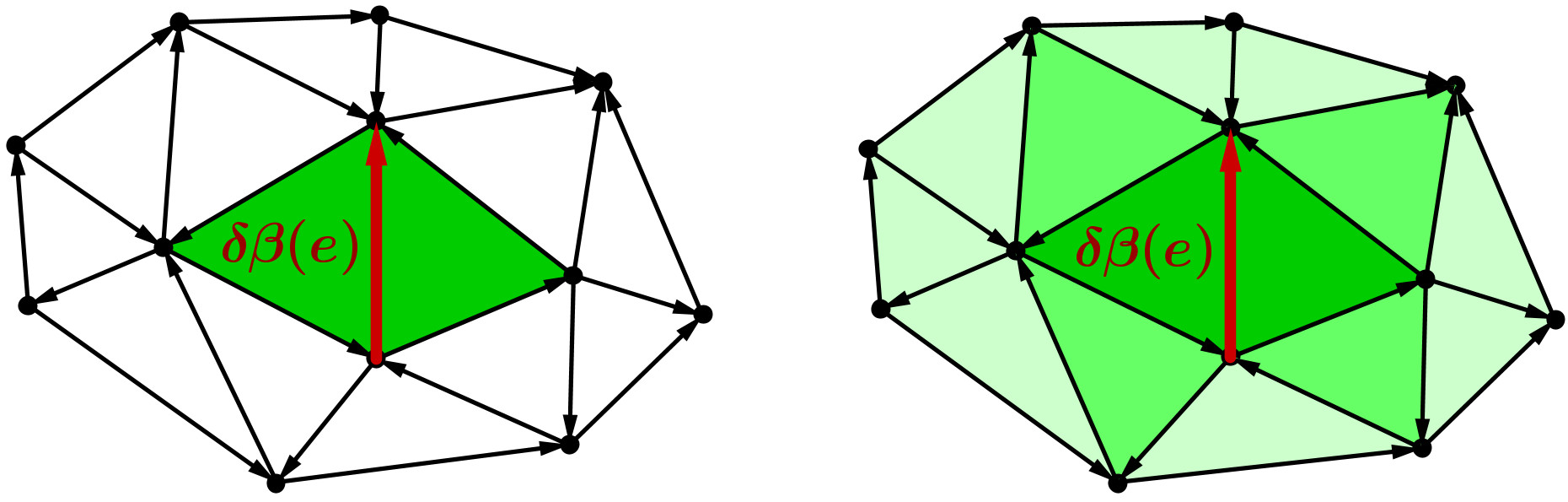}
\caption{Comparison of the support of the codifferential of 2--forms between the classical DEC (L) and our approach (R). The codifferential of a 2--form $\beta$ is a 1--form $\delta\beta$ located on edges. The value of $\delta\beta$ on the red edge $e$ is a linear combination of the values of $\beta$ on faces colored green.
The color intensity of faces reflects their weight of influence on $\delta\beta(e)$.
}\label{fig:codifferential2}
\end{figure}

In classical DEC the discrete codifferential operator on triangle meshes is defined using the diagonal approximation of the Hodge star operator. \citet{Alexa2011} in Section 3 hint at a codifferential of 1--forms on general polygonal meshes. The main difference between these and our codifferentials is in the support, see Figures \ref{fig:codifferential1} and \ref{fig:codifferential2}.

\subsection{The Laplace operator}
In exterior calculus, the Laplace operator is given by
$\Delta :=\delta d + d\delta,$
where $\delta$ is the codifferential and $d$ the exterior derivative. The Laplacian is defined in this way also in the classical DEC and we follow this convention. The classical Laplacian on 0--forms is also called the \textit{cotan Laplace operator}. Even though many different approaches lead to the cotan--formula, \citet{MacNeal} was the first to derive it.

\begin{figure}[htb]
\centering
\includegraphics[width=0.48\linewidth]{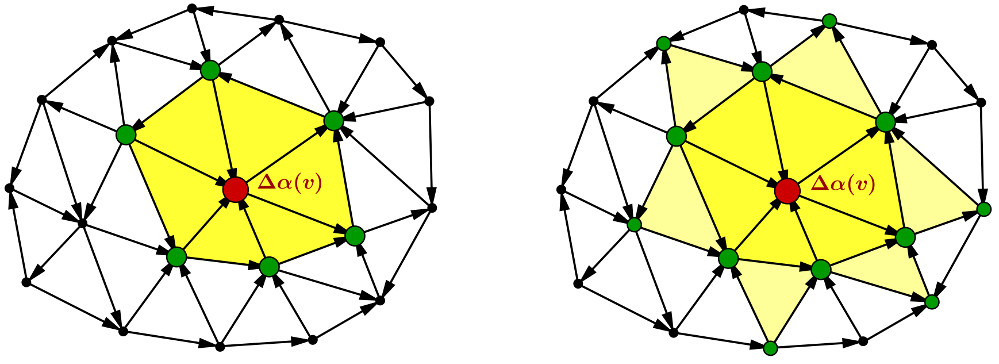}\hfill
\includegraphics[width=0.48\linewidth]{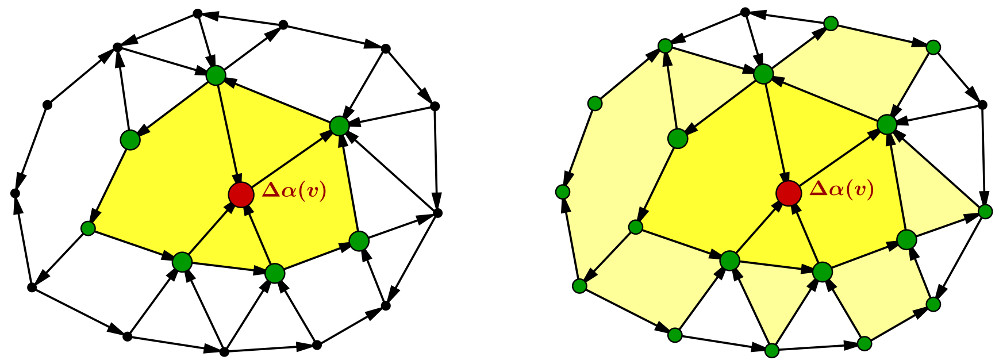}
\caption{Comparison of support of the Laplacian of 0--forms between the classical DEC (far left) and ours on triangle meshes (center left), Laplacian of \citet{Alexa2011} for their $\lambda=0$ (center right) and ours on polygonal meshes (far right). The Laplacian of a 0--form $\alpha$ is a 0--form $\Delta\alpha$ located again on vertices. The value of $\Delta\alpha$ on the red vertex $v$ is a linear combination of values of $\alpha$ on vertices colored green.
The support of our Laplacian is always larger, the point size reflects the weight of respective $\alpha$s on $\Delta\alpha(v)$. We also color yellow the faces whose vertices carry the $\alpha$s that enter as variables for $\Delta\alpha(v)$.
}\label{fig:triLaplacian0}
\end{figure}

Discrete Laplacians of 0--forms on general polygonal meshes were introduced by \citet{Alexa2011}. They derive Laplacians with a geometric and a combinatorial part that improves their mesh processing methods, see also Figure \ref{fig:LaplaceFlowKitty}.
By Theorem 2 therein, on triangle meshes their Laplacians reduce to the cotan--formula. In Figure \ref{fig:triLaplacian0} we compare the support of their purely geometric Laplacians to ours.

\subsection{The contraction operator and the Lie derivative}
The Lie derivative can be thought of as an extension of a directional derivative of a function to derivative of tensor fields (such as vector fields or differential forms) along a vector field. It is invariant under coordinate transformations, which makes it an appropriate version of a directional derivative on curved manifolds. It evaluates the change of a tensor field along the flow of a vector field and is widely used in mechanics.

Our discretization of Lie derivative of functions (0--forms) corresponds to the functional map framework of \citet{Azencot2013}, but now generalized to polygonal meshes. Our discrete Lie derivatives are thus linear operators on functions that produce new functions, with the property that the derivative of a constant function is 0.
Furthermore, both theirs and our Lie derivative of a vector field produces a vector field. However, whereas their framework is built for triangle meshes, we work with general polygonal meshes.

While maintaining the discrete exterior calculus framework, our work can also be interpreted as an extension of the Lie derivative of 1--forms presented by \citet{LieAdvection} from planar regular grids to surface polygonal meshes in space.


\section{Primal-to-Primal Operators}\label{chap:Contribution}

This section contains actual results of our research -- we present the theory and numerically evaluate the quality of our approximations by setting our results against analytical solutions. We also compare our methods to other DEC schemes.

Not using dual meshes simplifies the definition of several operators on polygonal meshes, which may be a difficult task otherwise. Moreover, it helps to maintain the compatibility of our operators since both the initial and the mapped discrete forms are located on primal elements. However, this approach also brings some drawbacks, we discuss them in this section.

\subsection{The Discrete Wedge Product}\label{subsec:CupProduct}
\begin{figure}[ht]
\centering\includegraphics[width=0.65\linewidth]{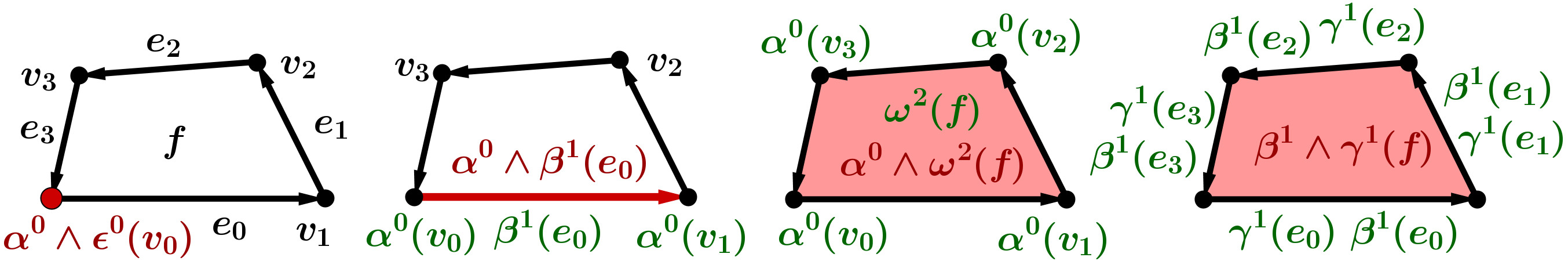}
\caption{The wedge product on a quadrilateral: the product of two 0--forms is a 0--form located on vertices (far left). The product of a 0--form with a 1--form is a 1--form located on edges (center left). The product of a 0--form with a 2--form is a 2--form located on faces (center right), and the product of two 1--forms is a 2--form located on faces (far right). }\label{fig:CupProductQuad}
\end{figure}

Just like the wedge product of differential forms, our discrete wedge product is a product of two discrete forms of arbitrary degree $k$ and $l$ that returns a form of degree $k+l$ located on primal $(k+l)$--dimensional cells
(see Figure \ref{fig:CupProductQuad}).

On triangle meshes our discrete wedge product is identical to the cup product given by \citet{Whitney} and on quadrilaterals it is equivalent to the cubical cup product of \citet{Rachel}. Further, the wedge product of differential forms satisfies the Leibniz product rule with exterior derivative and is skew--commutative. The discrete wedge product must satisfy these properties as well, we thus appropriately extend the discrete wedge product from triangles and quads to general polygons and derive the following formulas:
\begin{myDef}\label{def:PolygonalCup}
Let $S$ be a surface mesh (pseudomanifold) whose faces (2--cells) are polygons. The polygonal wedge product of two discrete forms $\alpha^k,\beta^l$ is a $(k+l)$--form $\alpha^k\wedge\beta^l$ defined on each $(k+l)$--cell $c^{k+l}\in S$.
Let $v$ be a vertex, $e=(v_i,v_j)$ an edge, and $f=(v_0,\dots,v_{p-1})$ a $p$--polygonal face with boundary edges $e_0,\dots,e_{p-1}$ having the same orientation as $f$. The \textbf{polygonal wedge product} is given for each degree by:
\begin{eqnarray*}
(\alpha^0\wedge\beta^0)(v)  &=&  \alpha(v)\beta(v),\\
(\alpha^0\wedge\beta^1)(e)  &=&  \frac{1}{2}(\alpha(v_i)+\alpha(v_{j}))\beta(e),\\
(\alpha^0\wedge\beta^2)(f) &=&\frac{1}{p}\bigg( \sum_{i=0}^{p-1}\alpha(v_i) \bigg) \beta(f),\\
(\alpha^1\wedge\beta^1)(f) &=& \sum_{a=1}^{\lfloor\frac{p-1}{2} \rfloor} \bigg( \frac{1}{2} - \frac{a}{p} \bigg)\sum_{i=0}^{p-1} \alpha(e_i)(\beta(e_{i+a}) - \beta(e_{i-a})), \text{where all indices are modulo }p.
\end{eqnarray*}
\end{myDef}

The polygonal wedge product is illustrated in Figures \ref{fig:CupTriangle11} and \ref{fig:CupProductQuad}.
It is a skew--commutative bilinear operation: 
$
\alpha^k \wedge\beta^l = (-1)^{kl} \beta^l\wedge\alpha^k,
$
matching its continuous analog, and it satisfies the Leibniz product rule with discrete exterior derivative:
$d(\alpha^k \wedge\beta^l) = d\alpha\wedge\beta + (-1)^{k}\alpha\wedge d\beta.$

The wedge product of three 0--forms is trivially associative (it reduces to multiplication of three scalars). Unfortunately for higher degree forms it is not associative in general, only if one of the 0--forms involved is constant. This is a common drawback of discrete wedge products, see e.g. \citep[Remark 7.1.4.]{Hirani}.

In matrix form, our polygonal wedge product reads:
\begin{align*}
 \alpha^0 \wedge \epsilon^0 &= \alpha^0 \odot \epsilon^0,\\
 \alpha^0 \wedge \beta^1 &= (\B \alpha^0) \odot \beta^1,\\
 \alpha^0 \wedge \omega^2 &= (\fv \alpha^0) \odot \omega^2,\\
 (\beta^1 \wedge \gamma^1)|_f &= (\beta^1|_f)^\top \R (\gamma^1|_f),
\end{align*}
where $\odot$ is the Hadamard (element--wise) product, $\beta|_f$ denotes the restriction of $\beta$ to a $p$--polygonal face $f$, and the matrices $\B\in\mathbb{R}^{|E|\times|V|}$, $\fv\in\mathbb{R}^{|F|\times|V|}$, and $\R\in\mathbb{R}^{p\times p}$ per $f$ read:
\begin{align}
\B[i,j] =& \left\{
      \begin{array}{ll}
	\frac{1}{2} & \mbox{if } v_j \prec e_i,\\
	0 & \mbox{otherwise.}
      \end{array}
      \right.  \label{eq:cup01} \\
\fv[i,j] =& \left\{
    \begin{array}{ll}
	\frac{1}{p_i} & \mbox{if } v_j \prec f_i, f_i\text{ is a } p_i\text{--gon}, \\
	0 & \mbox{otherwise.}
    \end{array}
    \right. \label{eq:cup02} \\
\R =& \sum_{a=1}^{\lfloor\frac{p-1}{2} \rfloor} \bigg( \frac{1}{2} - \frac{a}{p} \bigg) \R_a,
     \quad
\R_a[k,j] = \left\{
    \begin{array}{ll}
	1  & \mbox{if } e_j \text{ is $(k+a)$--th halfedge of } f, [f:e_j]=1, \\
	-1 & \mbox{if } e_j \text{ is $(k-a)$--th halfedge of } f, [f:e_j]=1, \\
	0 & \mbox{otherwise. }
    \end{array}
 \right. \label{eq:cup11}    
\end{align}


\subsubsection{Numerical evaluation}\label{subsub:NumericalEvaluationWedge}
We perform the numerical evaluation of our polygonal wedge product as an approximation to the continuous wedge product on a given mesh $S$ over a smooth surface in the following fashion:

\paragraph{1} We integrate each differential $l$--form over all $l$--dimensional cells of the mesh $S$ and thus define discrete forms $\alpha^0$, $\beta^1$, $\gamma^1$, and $\omega^2$:
 \[
  \alpha^0(v) = A(v),\;\; \beta^1(e) = \int\limits_e B,\quad \gamma^1(e)=\int\limits_e \Gamma,\quad
  \omega^2(f)= \int\limits_f\Omega,
 \]
where Greek capital letters denote the respective continuous differential forms. In practice, we integrate the continuous differential 2--form $\Omega$ over a set of triangles $(C,v_i,v_{i+1})$ that approximate the possibly non--planar face $f=(v_0,\dots,v_{p-1})$, where $C$ is the centroid of $f$, i.e.,
\[
 \omega^2(f)=\sum_{i=0}^{p-1} \int\limits_{(C, v_i,v_{i+1})}\Omega,\; \text{ where } C=\frac{1}{p}\sum_{i=0}^{p-1}v_i.
\]

\paragraph{2} Next we compute the polygonal wedge products 
 $(\alpha^0\wedge\beta^1)(e)$, $(\alpha^0\wedge\omega^2)(f)$, $(\beta^1 \wedge \gamma^1)(f)$ for each edge and face of the mesh using our formulas.

\paragraph{3} We also calculate analytical solutions  of the (continuous) wedge products and discretize (integrate) these solutions.

\paragraph{4} We then compute the $L^\infty$ and $L^2$ errors of our approximation. So let $\xi^k$ denote our solution (a discrete $k$--form) and $\Xi^k$ the respective discretized analytical solution, we compute:
\[
  L^2 \text{ error } = \Big(\xi^k - \Xi^k \Big)^\top M_k \Big(\xi^k - \Xi^k),\quad\quad
  L^\infty \text{ error } = ||\xi^k - \Xi^k ||_\infty = \max_{c^k}(|\xi^k(c^k) - \Xi^k(c^k)|),
\]
where $M_k$ are inner product matrices. Concretely, $M_2\in\mathbb{R}^{|F|^2}$, $M_0\in\mathbb{R}^{|V|^2}$ are diagonal matrices given by

\begin{equation}\label{eq:AWM0}
M_2[i,i]=\frac{1}{|f_i|},\quad  M_0[i,i]=\sum_{f_j\succ v_i}\frac{|f_j|}{p_j},
\end{equation}
and $M_1$ is the inner product of two 1--forms of \citet{Alexa2011}, i.e., for two 1--forms $\epsilon$ and $\lambda$, $M_1$ is defined in the sense that
\begin{equation}\label{eq:AWM1}
 \epsilon^\top M_1\lambda=\sum_f \epsilon|_f^\top M_f \lambda|_f,\quad
 M_f:= \frac{1}{|f|}B_f B_f^\top,
\end{equation}
where $\epsilon|_f$ again denotes the restriction of $\epsilon$ to a $p$--polygonal face $f$ and $B_f$ denotes a $p\times 3$ matrix with edge midpoint positions as rows (we take the centroid of each face as the center of coordinates per face).

\paragraph{5}
To evaluate the numerical convergence behavior, we refine the mesh over the given smooth underlying surface. The smooth surfaces used for tests are: unit sphere, torus azimuthally symmetric about the z-axis, and planar square. To create unstructured meshes, we randomly eliminate a given percentage of edges of an initially regular mesh.

We also use \textbf{jittering} to evaluate the influence of irregularity of a mesh on the experimental convergence. When jittering, we start with a regular mesh and displace each vertex in a random tangent direction to distance $r\cdot |e|$, where $|e|$ is the shortest edge length, and then project all thus displaced vertices on a given underlying smooth surface. 
That is why we use simple surfaces such as spheres and tori for our tests.

If not stated otherwise, all graphs use $\log_{10}$ scales on both the horizontal and vertical axes.

\begin{figure}[ht]
\centering
\includegraphics[width=0.25\linewidth]{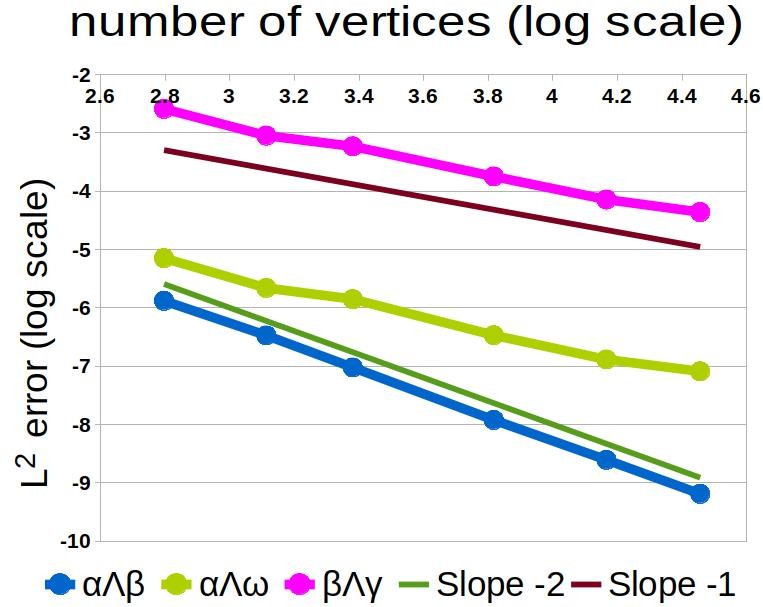}
\includegraphics[width=0.25\linewidth]{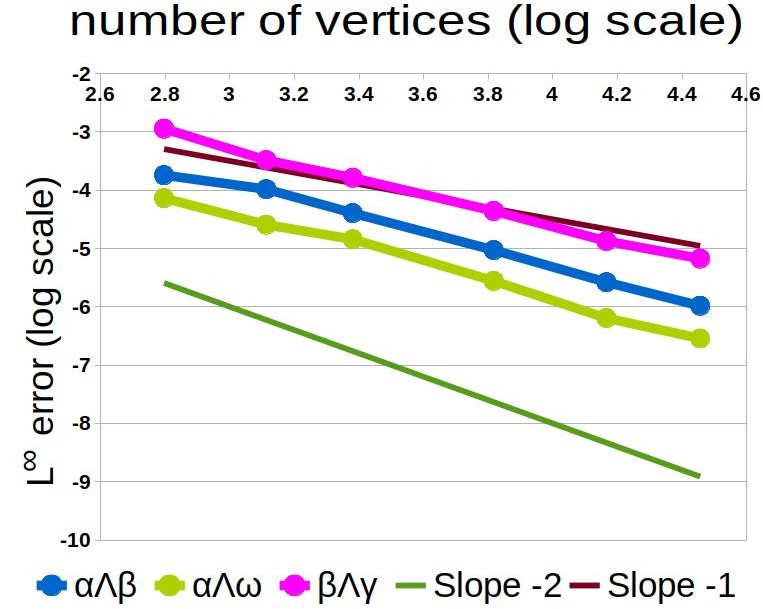}\hspace{0.02\linewidth}
\includegraphics[width=0.20\linewidth]{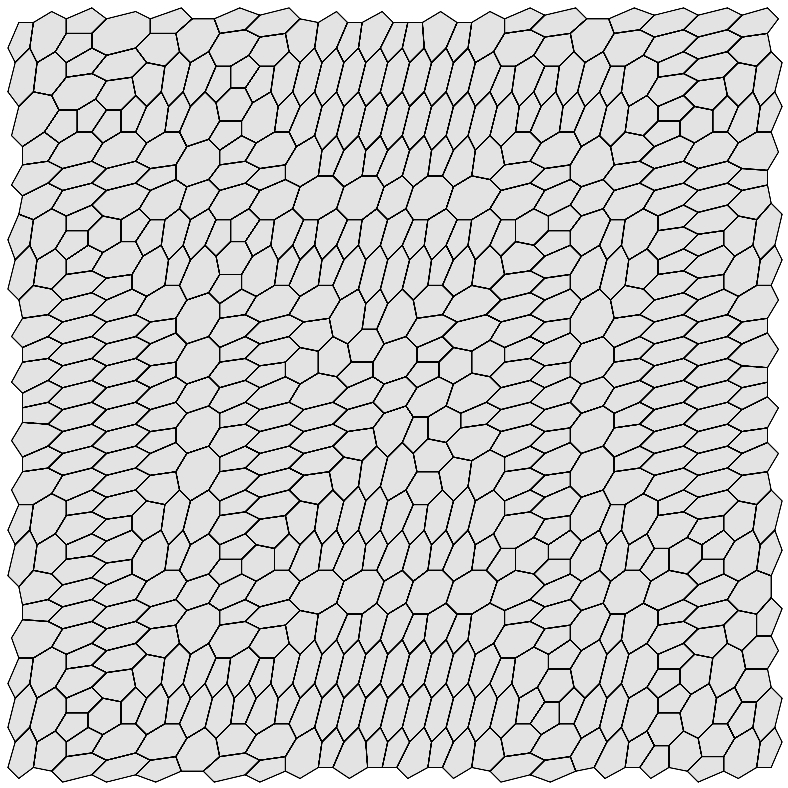}\hspace{0.01\linewidth}
\includegraphics[width=0.20\linewidth]{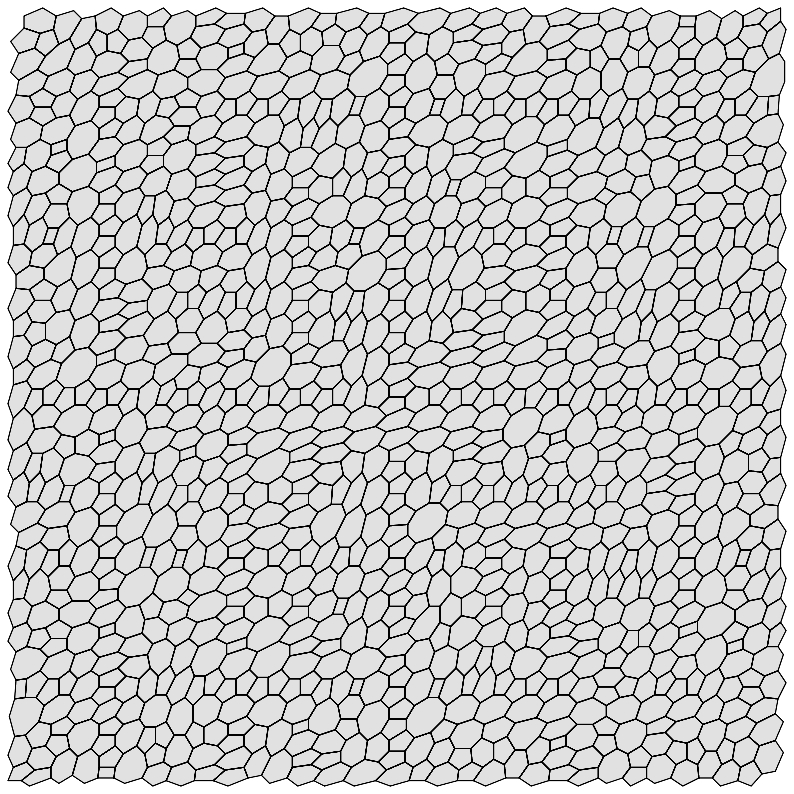}
\caption{Convergence of the wedge products on a set of unstructured polygonal meshes on a planar square to analytical solutions in $L^2$ norm (far left) and $L^\infty$ norm (center left). Both axes are in $\log_{10}$ scale.
The differential forms tested are trigonometric forms 
$\alpha^0 = \sin (x) \cos (y)+1$,
$\beta^1 = (\sin ^2(x)-1)dx  +(3\cos(x+2)+\sin (y))dy$,
$\gamma^1 = (\cos (x) \sin (y)+3) dx  + \cos (y)dy$,
$\omega^2 = (\sin (x y)+\cos (1)) dx \wedge dy$.
On the right are samples of tested meshes, both over a planar $[-1,1]^2$ square.}\label{fig:cupPolySquare}
\end{figure}

We have tested quadratic and trigonometric differential forms on flat and curved surface meshes (with non--planar faces) and our polygonal wedge products exhibit at least linear convergence to the respective analytical solutions, both in $L^2$ and $L^\infty$ norm. 
In Figure \ref{fig:cupPolySquare} we give an example.

\subsection{The Hodge Star Operator}\label{subsec:HodgeStar}

We define a discrete Hodge star operator as a homomorphism (linear operator) from the group of $k$--forms to $(2-k)$--forms. But since we do not employ any dual mesh and there is no isomorphism between the groups of $k$-- and $(2-k)$--dimensional cells, our Hodge star is not an isomorphism (invertible operator), unlike its continuous counterpart and diagonal approximations.

On the other hand, thanks to the dual forms being located on elements of our primal mesh, we can compute discrete wedge products of primal and dual forms and thus define a discrete inner product and discrete contraction operator later on.

Moreover, thanks to the Hodge star operating on primal meshes, we circumvent the ambiguity of defining dual meshes of unstructured general polygonal meshes. The idea of defining a Hodge star operator without using a dual mesh was borrowed from \citet{Rachel}, where the author suggests metric--independent Hodge star operators on simplicial and cubical complexes.

Our formulas for discrete Hodge star operators $\star$ are further motivated by the condition that the Hodge dual of constant discrete forms on planar surfaces is exact, hence $\star\mu=1$ and $\star 1=\mu$, where $\mu$ is the volume form on a given Riemannian manifold, for details see \citep[Section 6.5]{Abraham}.

\begin{figure}[htb]
\centering\includegraphics[width=0.35\linewidth]{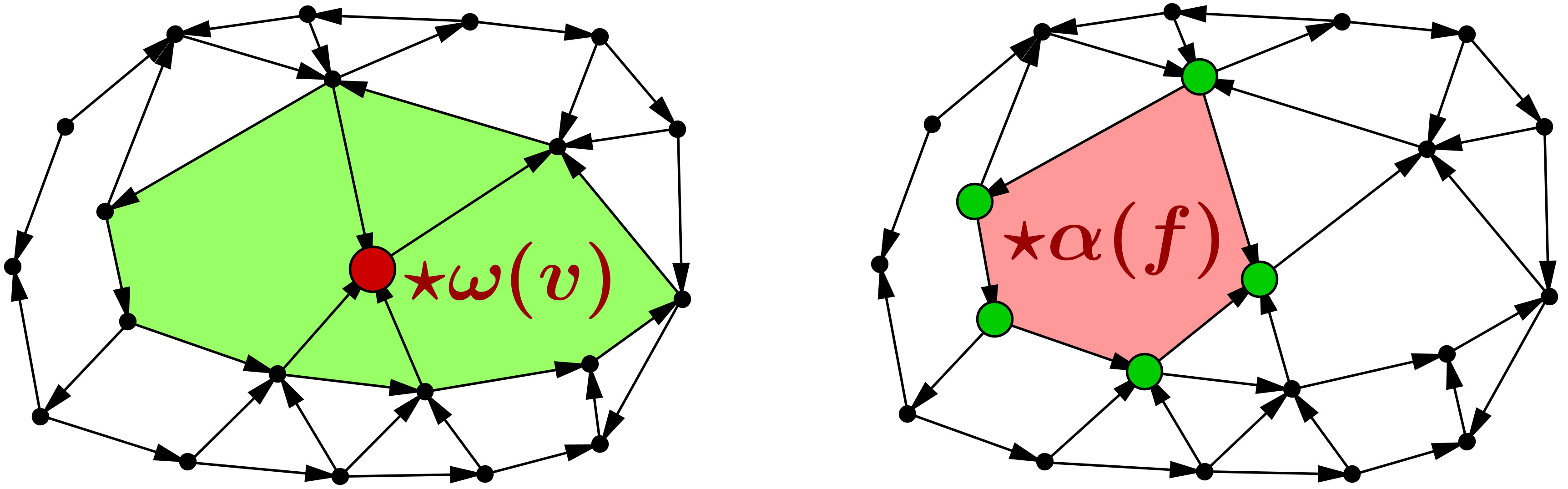}
\caption{On the left, the Hodge dual of a 2--form $\omega$ is a 0--form $\star\omega$, which value on a vertex $v$ (colored red) is a linear combination of values of $\omega$ on adjacent faces (colored green).
On the right, the Hodge dual of a 0--form $\alpha$ is a 2--form $\star\alpha$, the value of $\star\alpha$ on a face $f$ (colored red) is a linear combination of values of $\alpha$ on vertices (green) of that face.}\label{fig:Hodge02}
\end{figure}

\textbf{The Hodge star operator on 2--forms} takes in account the degree $p_i$ of $p_i$--polygonal faces $f_i$ and their vector areas $|f_i|$. If $\omega^2$ is a 2--form, then the 0--form $\star\omega$ on a vertex $v$ is given by
\small\begin{equation}
(\star_2\omega)(v)  = \frac{1}{\displaystyle\sum_{f_i\succ v}\frac{|f_i|}{p_i}} \cdot \sum_{f_i\succ v}\frac{\omega(f_i)}{p_i} ,\label{eq:DerivedHodge2}
\end{equation}\normalsize
i.e., it is a linear combination of values of $\omega$ on faces adjacent to $v$, see Figure \ref{fig:Hodge02} left.

\textbf{The Hodge star on an 1--form} $\beta^1$ is first defined per halfedges of a $p$--polygonal face $f$ as:
\begin{equation}
 \star_1\beta = \W_1\R^\top\beta,\label{eq:DerivedHodge1}
\end{equation}
where $\R$ is the matrix defined in (\ref{eq:cup11}) and $\W_1$ is a symmetric $p\times p$ matrix given by:
\small\[
\W_1[i,j] = \frac{\langle e_i,e_j \rangle}{|f|},
\]\normalsize
for $e_k$ the halfedges incident to and having the same orientation as the face $f$, where $\langle .,.\rangle$ denotes the Euclidean dot product. 

\begin{figure}[htb]
\centering\includegraphics[width=0.16\linewidth]{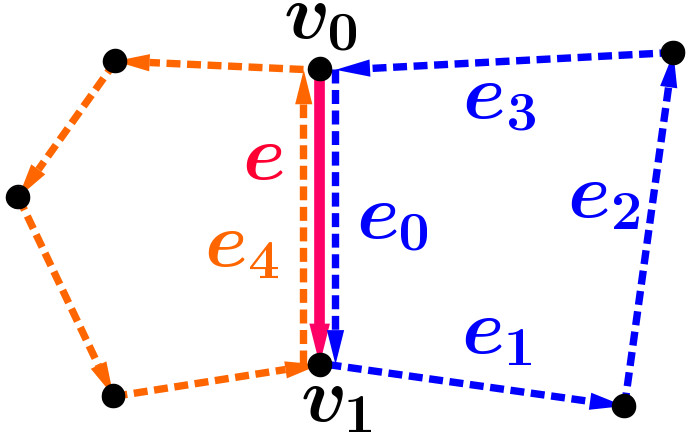}
\caption{Let $\beta\in C^1$ and $e=(v_0,v_1)$ be the edge with $e_0,e_4$ as the corresponding halfedges, then $\star\beta(e) = \frac{\star\beta(e_0)-\star\beta(e_4)}{2}$, where $\star\beta(e_4)$ is a linear combination of values of $\beta$ on dashed orange halfedges and $\star\beta(e_0)$ is a linear combination of values of $\beta$ on dashed blue halfedges, concretely
\small
$
\star\beta(e_0) = \frac{1}{4|f_0|}\Big((\langle e_0,e_1\rangle-\langle e_0,e_3\rangle)(\beta(e_0)-\beta(e_2))
+ (\langle e_0,e_0\rangle-\langle e_0,e_2\rangle)(\beta(e_3)-\beta(e_1)\Big).
$
\normalsize
}\label{fig:Hodge1}
\end{figure}

If an edge $e$ is not on boundary, it has two adjacent faces and two halfedges, thus we compute the values of $\star\beta$ on corresponding halfedges, sum their values with appropriate orientation sign and divide the result by 2, see an example in Figure \ref{fig:Hodge1}.

\textbf{The Hodge dual of a 0--form} $\alpha$ is a 2--form $\star\alpha$ defined per a $p$--polygonal face $f$ by:
\begin{equation}
 (\star_0\alpha) (f)  =  \frac{|f|}{p}\sum_{v_i\succ f} \alpha(v_i),\label{eq:DerivedHodge0}
\end{equation}
and it is simply the arithmetic mean of the values of $\alpha$ on vertices of the given face $f$ multiplied by the vector area $|f|$.

In matrix form, the discrete Hodge star operators read
\begin{align*}
 \star_0 &= \W_F\fv ,\\
 \star_1 &= \A\W_1\R^\top ,\\
 \star_2 &= \W_V\fv^\top,
\end{align*}
where $\fv$ and $\R$ are defined in equations (\ref{eq:cup02}) and (\ref{eq:cup11}), resp., and
$\W_F\in\mathbb{R}^{|F|\times|F|}$, $\W_V\in\mathbb{R}^{|V|\times|V|}$, $\A\in\mathbb{R}^{|E|\times|E|}$ are given by

\[
 \W_F[i,i] =|f_i|,\quad \W_V[i,i] =\frac{1}{\sum\limits_{f_k\succ v_i}\frac{|f_k|}{p_k}},\quad
 \A[i,j] = \left\{
    \begin{array}{cl}
	1 		& \mbox{if } i=j,\; e_i \text{ is on boundary}, \\
	\frac{1}{2} 	& \mbox{if } i=j,\; e_i \text{ is not on boundary}, \\
	-\frac{1}{2}  	& \mbox{if } e_i = -e_j, \\
	0 & \mbox{otherwise.}
    \end{array}
 \right.
\]

\begin{figure}[htb]
\centering
\includegraphics[width=0.25\linewidth]{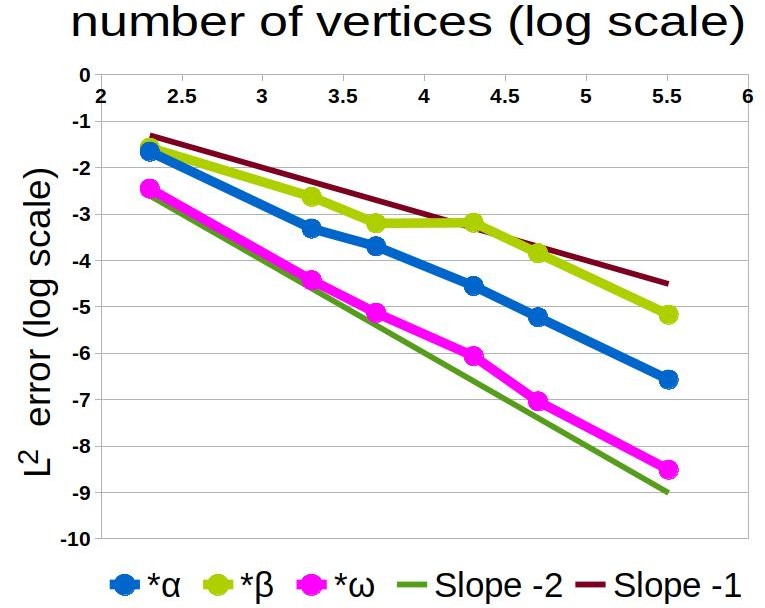}
\includegraphics[width=0.25\linewidth]{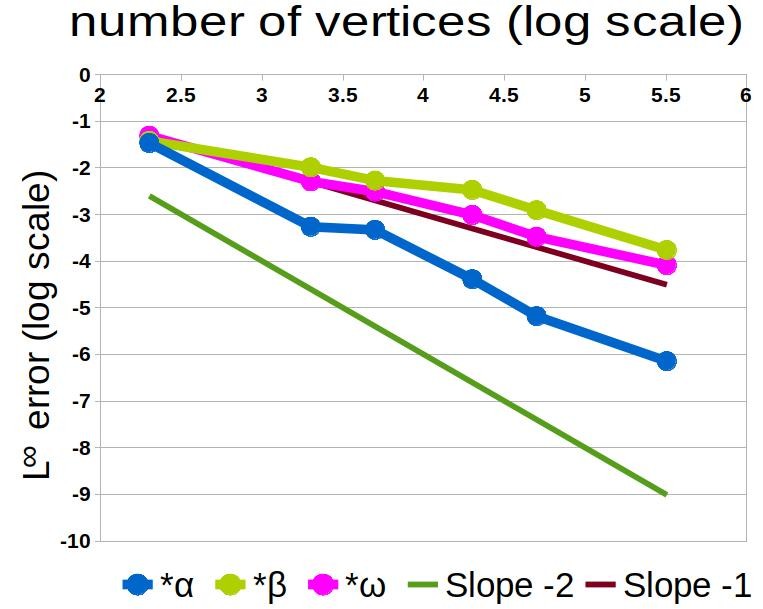}
\includegraphics[width=0.24\linewidth]{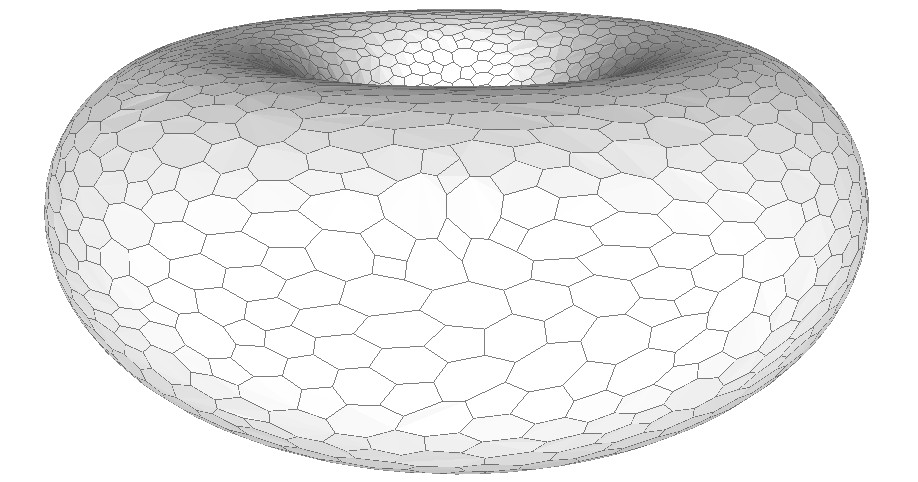}
\includegraphics[width=0.24\linewidth]{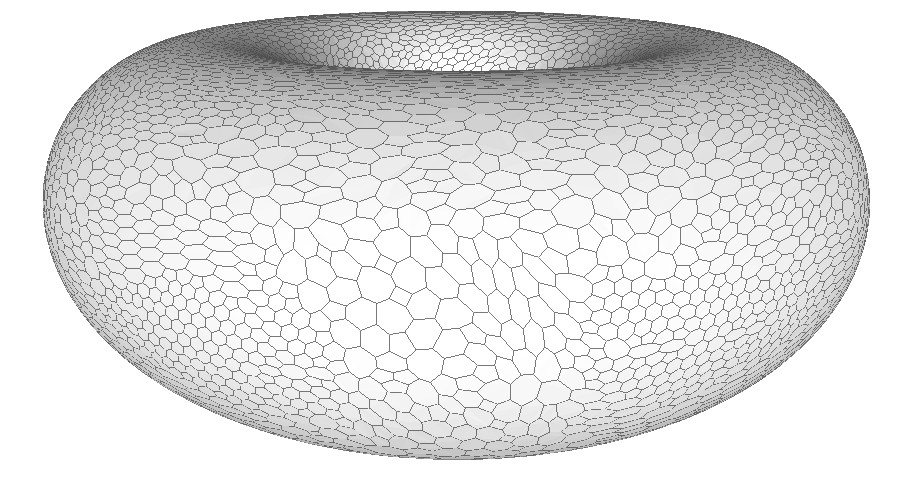}
\caption{The approximation errors of the discrete Hodge star on a set of irregular polygonal meshes on the torus. We have chosen $\alpha^0 =x^2 + y^2$, $\beta^1 = \mathbb{X}^\flat$, where $\mathbb{X}=(-y,x,0)$ is a tangent vector field, and $\omega^2=\mu$ is the area element on the torus. Thus $\star\mu = 1$, $\star\alpha = (x^2+y^2)\mu$, and $\star\beta = \mathbb{Y}^\flat$, where $\mathbb{Y}=2(-xz,-yz,x^2+y^2 -\sqrt{x^2+y^2})$ is a tangent vector field orthogonal to $\mathbb{X}$. On the right are two examples of meshes on the torus with 5k vertices and 20k vertices.}\label{fig:HodgeTorus}
\end{figure}

Although our Hodge star matrices are not diagonal, they are highly sparse and thus computationally efficient.
We have performed numerical tests on linear, quadratic, and trigonometric forms on planar and curved meshes and they exhibit the same at least linear convergence rate. We give an example in Figure \ref{fig:HodgeTorus}.

\subsection{The Hodge Inner Product}\label{subsec:HodgeInnerProduct}
The $L^2$--Hodge inner product of differential forms $\Gamma^k, \Omega^k$ on a Riemannian manifold $M$ is defined as: 
\[(\Gamma^k,\Omega^k) := \int_M\Gamma \wedge\star\Omega.\]
We define a \textbf{discrete $L^2$--Hodge inner product} of two discrete forms $\alpha^k, \beta^k$ on a mesh $S$ by:
\[
(\alpha^k,\beta^k) := \sum_{f\in S} \big(\alpha\wedge\star\beta\big)(f) =\alpha^\top\M_k \beta ,\; k=0,1,2,
\]
where $\M_k$ are the discrete Hodge inner product matrices that read:
\begin{align*}
 \M_0 =& \fv^\top\W_F\fv ,\\
 \M_1 =& \R\A\W_1\R^\top ,\\
 \M_2 =& \fv\W_V\fv^\top.
\end{align*}

It can be shown that our inner product of 1--forms restricted to a single face $f$ is identical to the one of \citep[Lemma 3]{Alexa2011}:
$
\R\A\W_1\R^\top|_f = M_f, 
$
where $M_f$ is defined as in equation (\ref{eq:AWM1}).
However, if a given mesh $S$ is not just a single face, they differ in general, i.e., for 1--forms $\beta^1$, $\gamma^1$:

\[
 \beta^\top\M_1\gamma = \beta^\top\R\A\W_1\R^\top\gamma
 \;\neq\; \beta^\top\R\W_1\R^\top\gamma = \sum_{f\in S} \beta^\top M_f \gamma.
\]

To numerically evaluate our inner products, we calculate our discrete $L^2$--Hodge norms of forms $\alpha^0, \beta^1, \omega^2$ over a mesh $S$ and compare them to their respective analytical $L^2$ norms. That is, if $\Gamma^k$ is a differential $k$--form and $\gamma^k$ the corresponding discrete $k$--form, we compute the error of approximation as:
\small\[
\int\displaylimits_S\Gamma\wedge\star\Gamma - \sum_{f\in S}\gamma\wedge\star\gamma = 
\int\displaylimits_S\Gamma\wedge\star\Gamma - \gamma^\top\M_k\gamma.
\]
\normalsize

\begin{figure}[htb]
\centering
\includegraphics[width=0.25\linewidth]{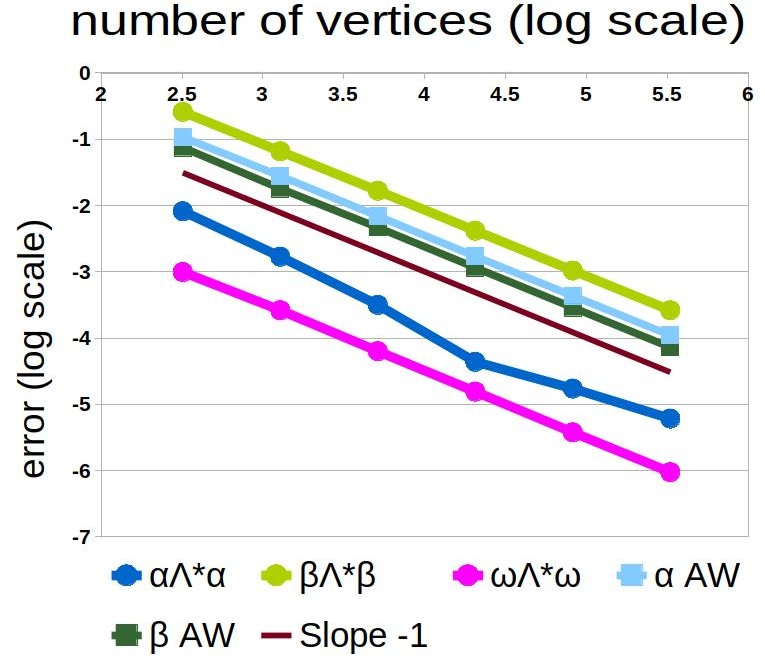}
\includegraphics[width=0.25\linewidth]{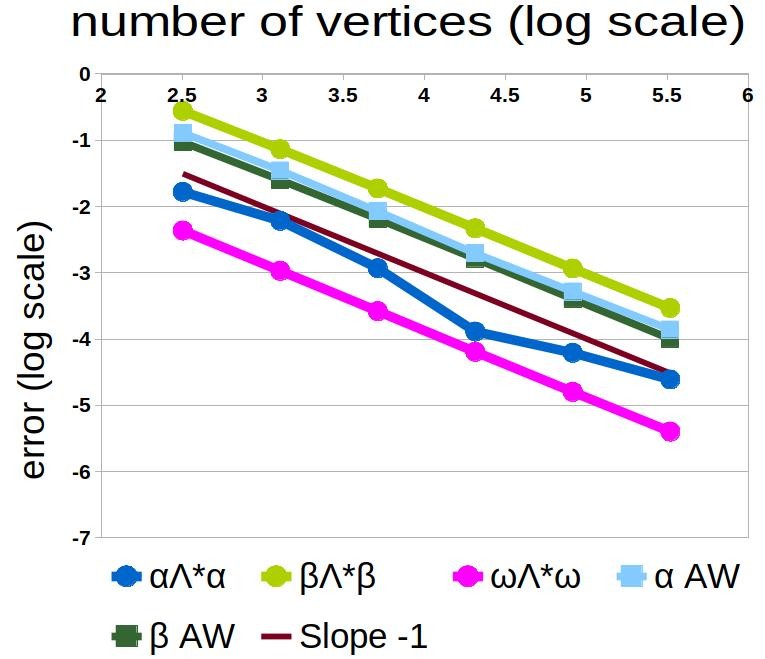} 
\includegraphics[width=0.23\linewidth]{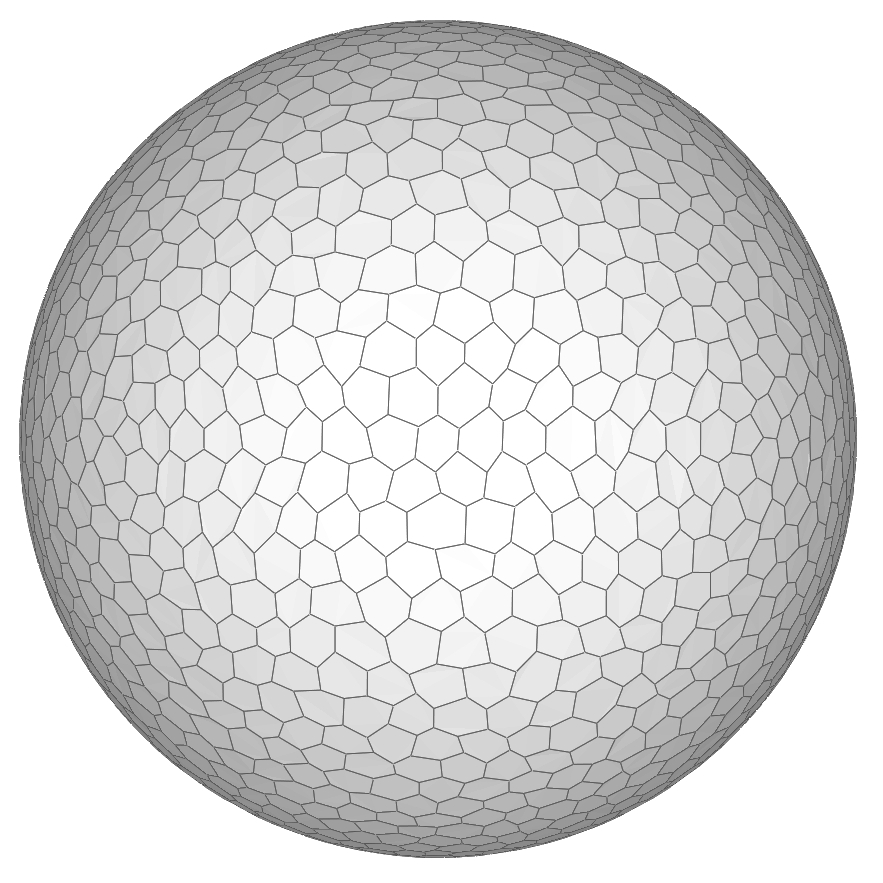}
\includegraphics[width=0.23\linewidth]{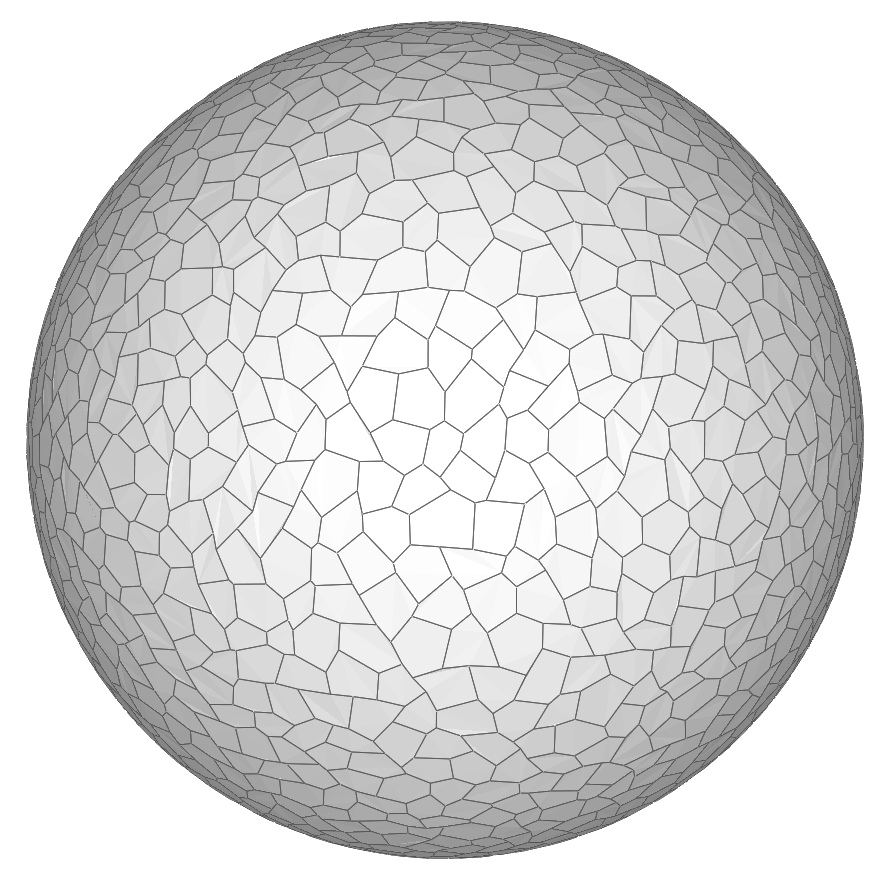}
\caption{The influence of irregularity on experimental convergence of discrete Hodge inner products to analytically computed solutions. Concretely, on a set of jittered meshes with $r=0.2$ (far left) and $r=0.4$ (center left). We show a sample mesh with 5k vertices for $r=0.2$ (center right) and $r=0.4$ (far right).
We set $\alpha^0=x^2+y^2$, $\beta^1=-xzdx -yzdy + (x^2+y^2)dz$, $\omega^2=xdy\wedge dz + y dz\wedge dx + z dx\wedge dy $. Here $\alpha\wedge\star\alpha$ denotes our Hodge inner product on $\alpha$, and similarly for $\beta$ and $\omega$. $\alpha\mathsf{AW}$ denotes the inner product of 0--forms and $\beta\mathsf{AW}$ the product of 1--forms of \citet{Alexa2011}.
}\label{fig:HodgeInnerProductSphere}
\end{figure}

An example of numerical evaluation of our $L^2$--Hodge inner products and numerical evaluation of inner products $M_0$ and $M_1$ of \citet{Alexa2011}, see also the equations (\ref{eq:AWM0} -- \ref{eq:AWM1}), is given in Figure \ref{fig:HodgeInnerProductSphere}.
The experimental convergence rate of our discrete $L^2$ norms is at least linear on all tested forms on compact manifolds with or without boundary.

\subsection{The Contraction Operator}\label{subsec:ContractionOperator}
The contraction operator $\I_X$, also called the interior product, is the map that sends a $k$--form $\omega$ to a $(k-1)$--form $\I_X\omega$ such that
$
 (\I_X\omega)(X_1,\dots,X_{k-1}) = \omega(X,X_1,\dots,X_{k-1})
$
for any vector fields $X_1,\dots,X_{k-1}$. The following property holds \citep[Lemma 8.2.1]{Hirani}:
\begin{myLemma}
Let $M$ be a Riemannian $n$--manifold, $X\in \mathfrak{X}(M)$ a vector field, then for the contraction of a differential $k$--form $\alpha$ with a vector field $X$ holds:
\[
 i_X\alpha = (-1)^{k(n-k)}\star(\star\alpha\wedge X^\flat),
\]
where $\flat: \mathfrak{X}(M)\rightarrow\Omega(M)$ is the flat operator.
\end{myLemma}

Since we already have discrete wedge and Hodge star operators that are compatible with each other, we can employ the lemma to define our \textbf{discrete contraction operator} $\I_X:C^k(S)\rightarrow C^{k-1}(S)$ on a polygonal mesh $S$ by:
\small\begin{equation}\label{eq:contraction}
\I_X\alpha = (-1)^{k(2-k)}\star(\star\alpha\wedge X^\flat),\;\;\alpha\in C^k(S),\; k = 1,2,
\end{equation}\normalsize
where the \textbf{discrete flat operator} on a vector field $X$ is given by discretizing its value over all edges of $S$. Let $e=(v_0,v_1)$ be an edge of $S$, then $e=e(t)=v_0+(v_1-v_0)t, t\in[0,1]$, and we set:
\begin{equation}\label{eq:flatOperator}
X^\flat(e) = \int\limits_e \langle e',X \rangle = \int\limits_{0}^1 \langle e'(t),X(e(t)) \rangle dt.
\end{equation}

Thus the discrete contraction operator is a linear operator that maps $k$--forms located on $k$--dimensional primal cells to $(k-1)$--forms located on $(k-1)$--dimensional primal cells.

\begin{figure}[htb]
\centering
\includegraphics[width=0.24\linewidth]{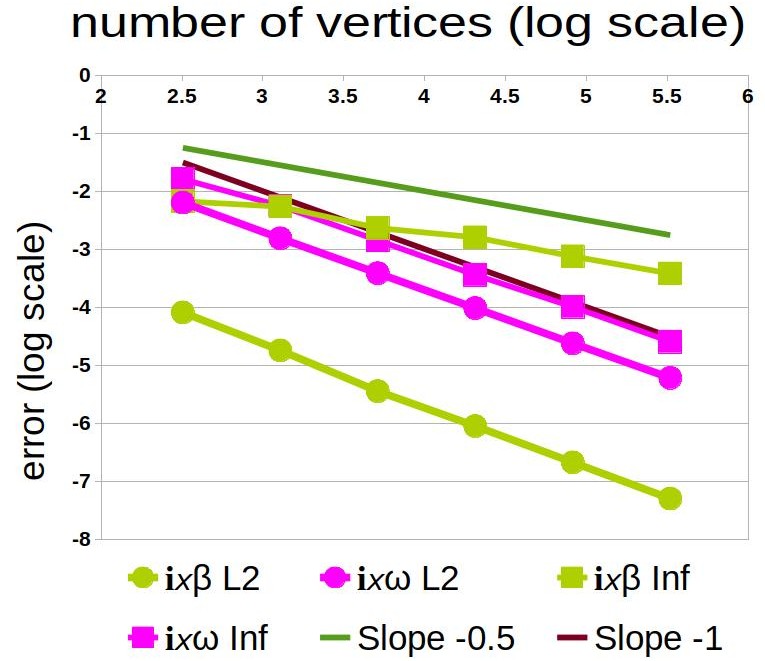}
\includegraphics[width=0.24\linewidth]{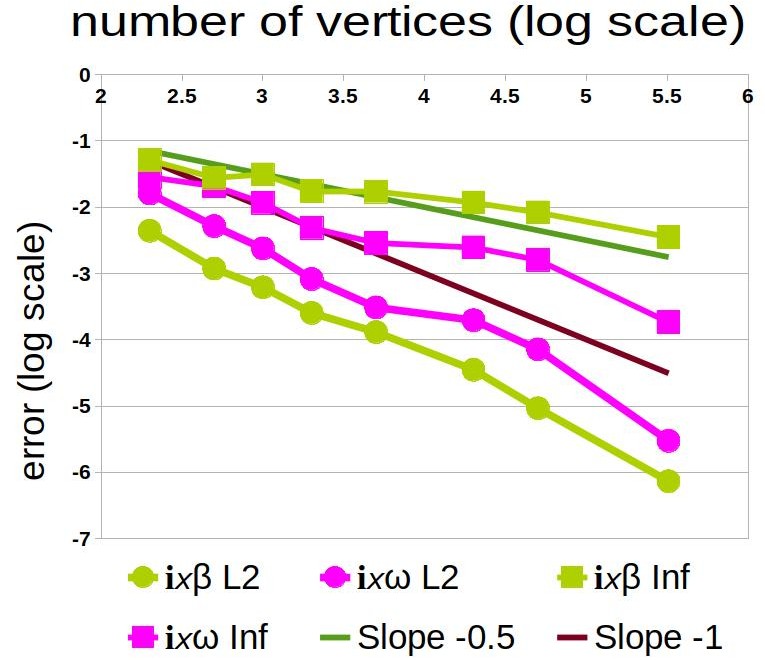}
\caption{The contraction operator on a unit sphere (L) and a torus (R). For the sphere we have used the same set of jittered meshes with $r=0.4$ as in Figure \ref{fig:HodgeInnerProductSphere}, and contracted the same forms as therein with respect to vector field $X = (-y,x,0)$. For the torus (R) we have contracted the differential forms of Figure \ref{fig:HodgeTorus} with respect to
vector field $X = 2(-xz, -yz, x^2 + y^2 -\sqrt{x^2 + y^2})$ on the same set of meshes as therein.
$\I_X\beta\;\mathsf{L2}$ denotes the $L^2$ error approximation of the contraction operator on the 1--form $\beta$, whereas $\I_X\beta\;\mathsf{Inf}$ denotes the $L^\infty$ error approximation on $\beta$, and similarly for the 2--form $\omega$.}\label{fig:contraction3D}
\end{figure}

Our discrete contraction of differential 2--forms with respect to different vector fields exhibit linear convergence to the analytically computed solutions, both in $L^\infty$ and $L^2$ norms. On 1--forms, the errors of approximation decrease linearly in $L^2$ and with slope $0.5$ in $L^\infty$ norm, see two examples in Figure \ref{fig:contraction3D}.

\subsection{The Lie Derivative}\label{subsec:LieDerivative}
We define the \textbf{discrete Lie derivative} $\Le_X:C^k(S)\rightarrow C^{k}(S)$ using Cartan's magic formula:
\begin{equation}\label{eq:discreteLie}
 \Le_X \alpha= \I_X d\alpha + d\I_X\alpha, \;\alpha\in C^k(S),\; k = 0,1,2.
\end{equation}
Unfortunately, the Leibniz product rule of our contraction operator and Lie derivative with discrete exterior derivative is not satisfied in general. Concretely
\begin{eqnarray*}
 \I_X(\alpha^k\wedge\beta^l)&=&(\I_X\alpha^k)\wedge\beta^l + (-1)^k \alpha^k\wedge(\I_X\beta^l) ,\\ \Le_X(\alpha^k\wedge\beta^l)&=&(\Le_X\alpha^k)\wedge\beta^l + \alpha^k\wedge(\Le_X\beta^l),
\end{eqnarray*}
holds only if $\alpha$ or $\beta$ is a closed 0--form.
Already \citet{Hirani} noticed that the Leibniz rule for Lie derivative might not hold due to the discrete wedge product not being associative in general. We confirm the observation of \citet{DEC2005} that the Leibniz rule may be satisfied only for closed forms.

\begin{figure}[ht]
\centering
\includegraphics[width=0.24\linewidth]{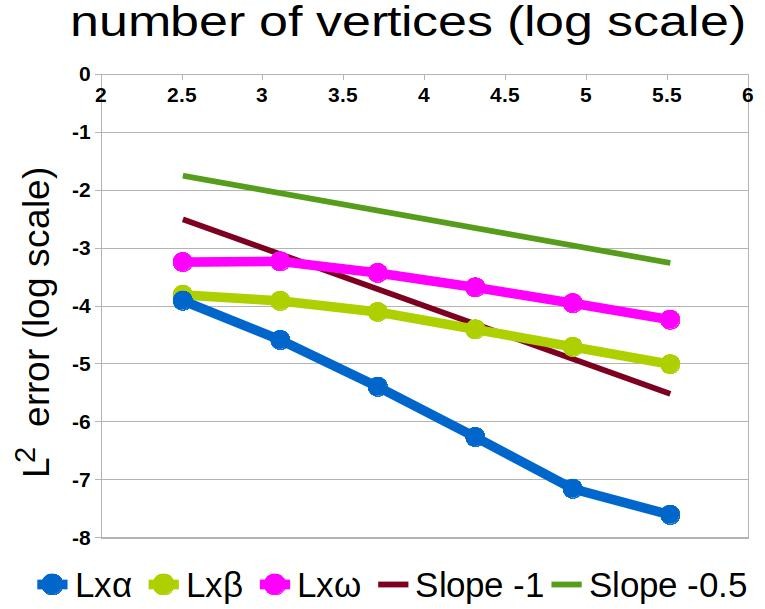}
\includegraphics[width=0.24\linewidth]{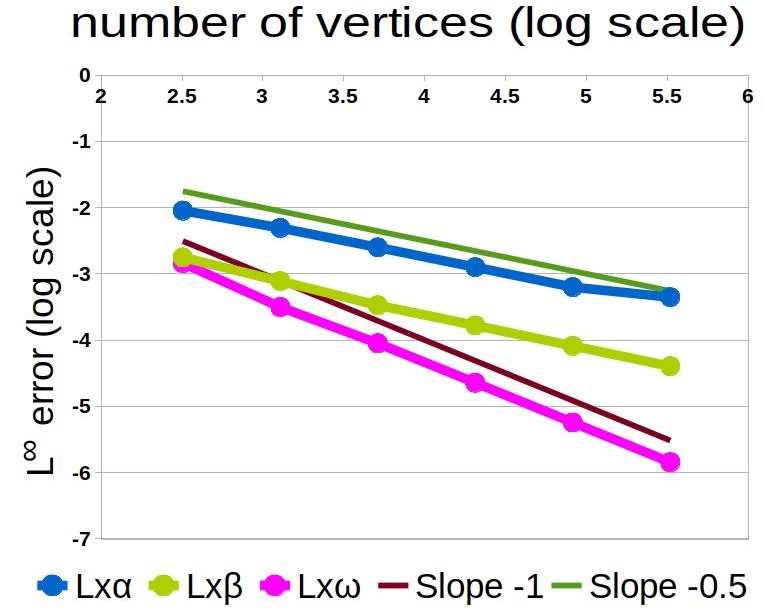}\hfill  
\includegraphics[width=0.24\linewidth]{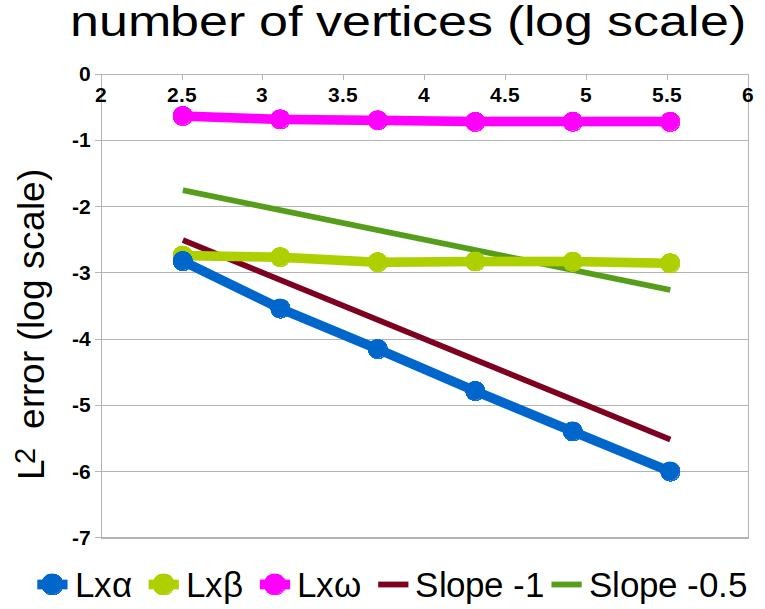}
\includegraphics[width=0.24\linewidth]{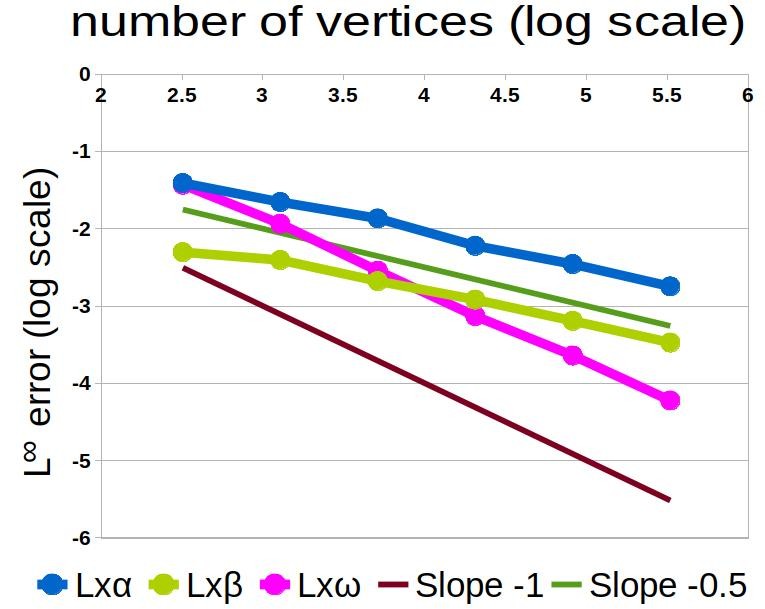}
\caption{The influence of jittering on experimental convergence of discrete Lie derivatives on regular polygonal meshes on a sphere (far and center left) and jittered meshes with vertex displacement by $0.4\times$ shortest edge length (center and far right). We use the same set of jittered meshes as in Figure \ref{fig:HodgeInnerProductSphere}. We employ the same forms and a vector field as in Figure \ref{fig:contraction3D}.
}\label{fig:LieDerivativeSphere}
\end{figure}

The Lie derivatives exhibit converging behavior on all tested forms on regular meshes, planar and non--planar. However, the $L^2$ error of approximation of Lie derivatives of 1-- and 2--forms on irregular meshes stays rather constant, see an example on a set of regular versus jittered meshes on a unit sphere in Figure \ref{fig:LieDerivativeSphere}. In this figure we can see that the $L^2$ error of the Lie derivative of a 1--form $\beta$ and a 2--form $\omega$ on regular meshes decreases with slope $-0.5$, whereas on very irregular meshes it stays constant.

Although our discrete Lie derivative of 1--forms on irregular meshes does not converge, in general, to analytically computed solutions, it can still be employed for Lie advection (Section \ref{subsec:LieAdvection}) of vector fields on irregular meshes and produce visually satisfying results, as we demonstrate in Figure \ref{fig:LieAdvectionTorus}.

\begin{figure}[htb]\centering
\includegraphics[width=0.15\linewidth]{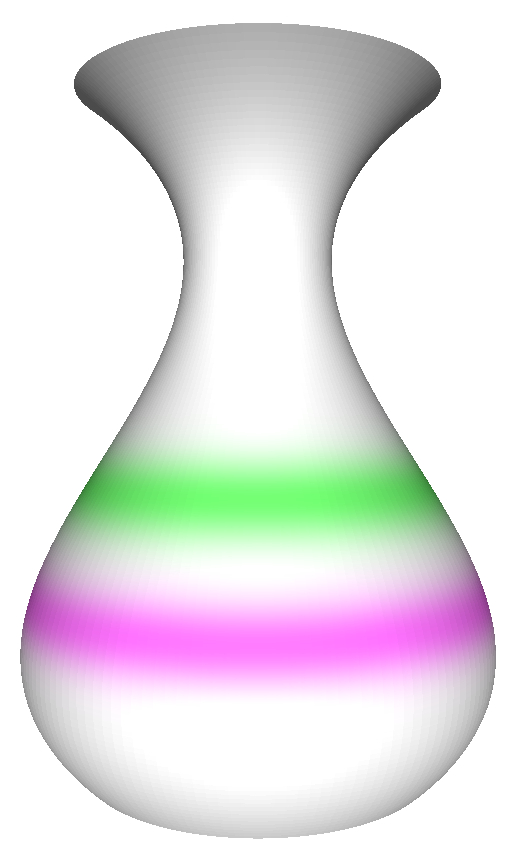}
\includegraphics[width=0.15\linewidth]{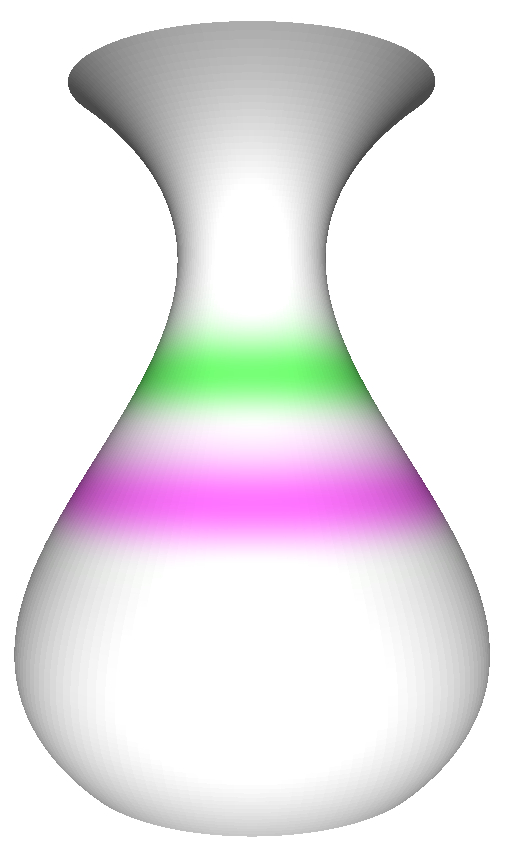}
\includegraphics[width=0.15\linewidth]{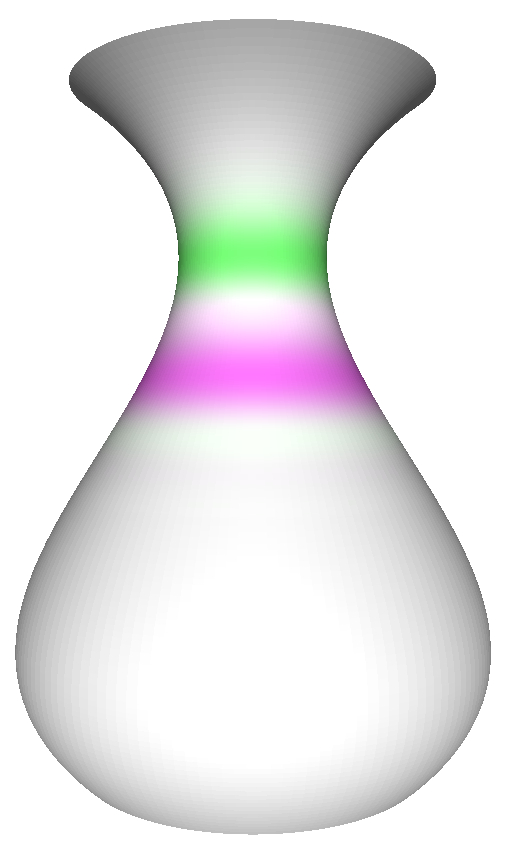}
\includegraphics[width=0.15\linewidth]{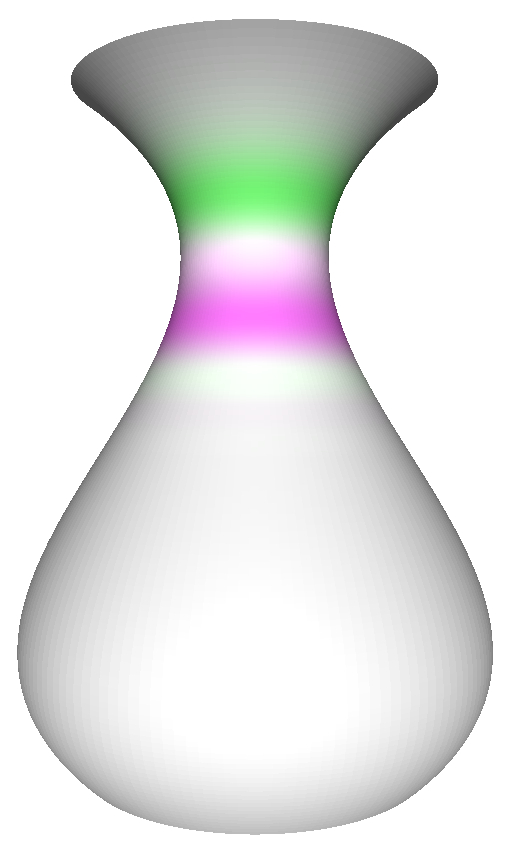}
\caption{Lie derivative used for advection of a color function (far left) encoded as a 0--form $\beta$ and advected using equation (\ref{eq:LieAdvection}). The advected function after 2000 (center left), 4000 (center right), and 5000 iterations (right) with time step $10^{-2}$. }\label{fig:LieVaso}
\end{figure}

\subsection{The Codifferential Operator}\label{subsec:Codifferential}
Just like on Riemannian $n$--manifolds, we define our \textbf{discrete codifferential operator} $\delta$ as
\[
 \delta_k\beta^k=(-1)^{n(k-1)+1}\star d\star\beta,\; \text{ $\beta$ a discrete $k$-form}.
\]
Thus in matrix form our codifferential operators read:
\begin{align*}
 \delta_1 =& - \W_V\fv^\top\, d_1\, \A\W_1\R^\top,\\
 \delta_2 =& - \A\W_1\R^\top\, d_0\, \W_V\fv^\top.
\end{align*}

If $M$ is a compact manifold without boundary or if $\alpha$ or $\star\beta$ has zero boundary values, then the codifferential is the adjoint operator of the exterior derivative with respect to the $L^2$--Hodge inner product:
$
 (d\alpha,\beta)=(\alpha,\delta\beta) \quad \forall \alpha\in\Omega^{k-1}(M),\beta\in\Omega^{k}(M).
$
\citet{Alexa2011} use this equation to derive their discrete codifferential operator on 1--forms, that reads
\[
 \delta_1 = M_0^{-1} d_0^\top M_1,
\]
where $M_0$ and $M_1$ are as in equations (\ref{eq:AWM0}--\ref{eq:AWM1}). This codifferential reduces to the classical codifferential, e.g., \citep{DDFormsForCM}, in the case of a pure triangle mesh.

\begin{figure}[htb]
\centering
\includegraphics[width=0.25\linewidth]{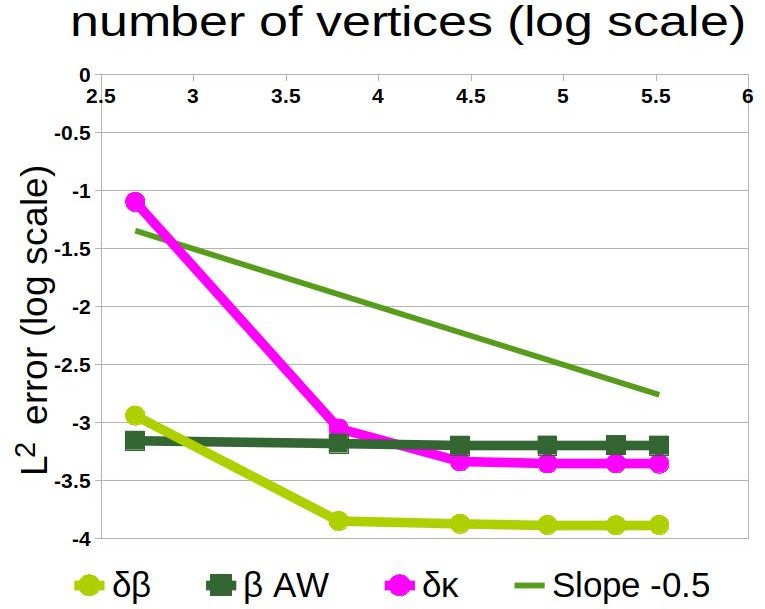} 
\includegraphics[width=0.25\linewidth]{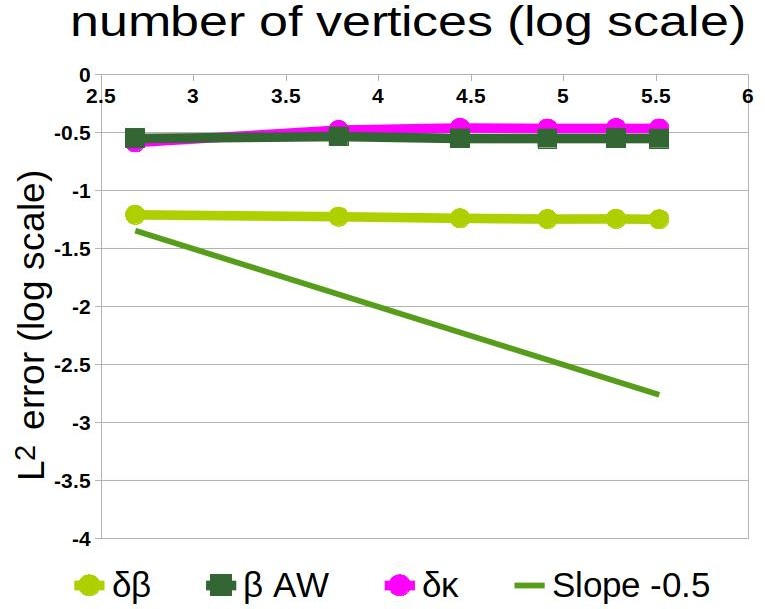}\hspace{0.02\linewidth}
\includegraphics[width=0.18\linewidth]{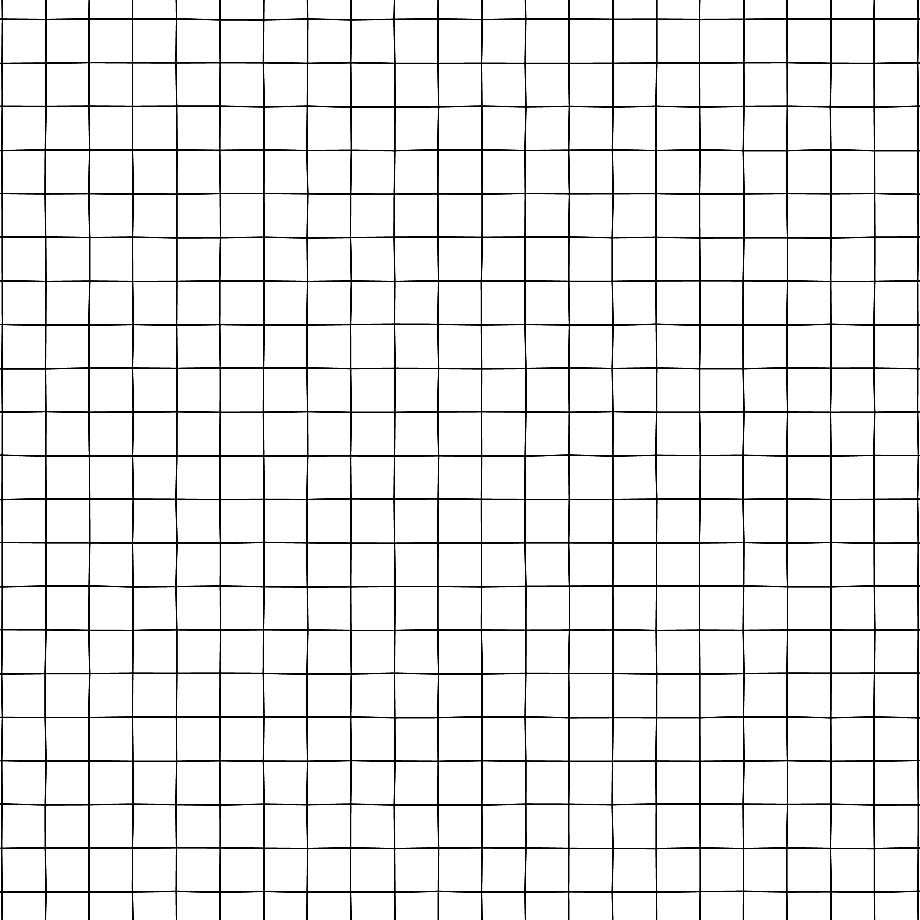}\hspace{0.01\linewidth}
\includegraphics[width=0.18\linewidth]{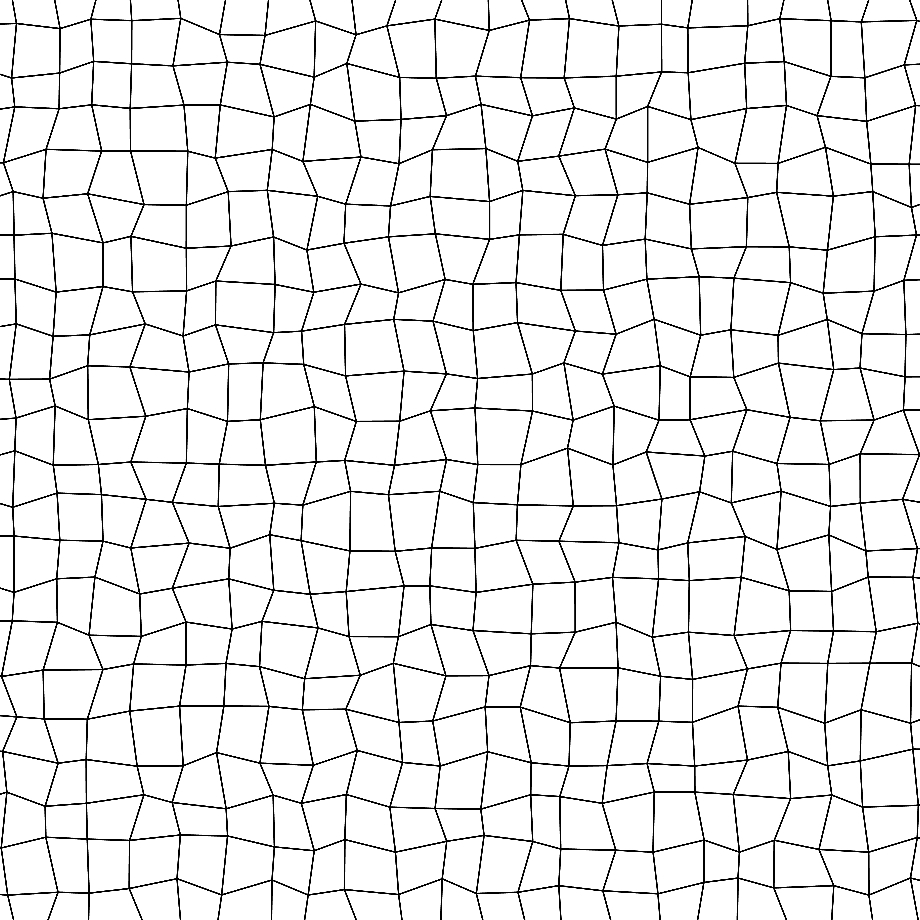}
\caption{The influence of jittering on experimental convergence of codifferentials of 
$\beta^1 =(\sin(2x) + \cos(\frac{y}{2}))dx + (3\sin(x) - \cos(y))dy$
and $\kappa^2  = (\sin(\frac{x+1}{4}) + \cos(1-\frac{y}{3}))dx\wedge dy$ on a set of planar quadrilateral jittered meshes with $r=0.01$ (far left) and $r=0.2$ (center left). $\delta\beta$ denotes the approximation error of our codifferential of $\beta$, and similarly for $\delta\kappa$. $\beta\mathsf{AW}$ stands for the $L^2$ error of the codifferential of \citet{Alexa2011}. On the center and far right are samples of such jittered meshes.
}\label{fig:codifferentialGrid}
\end{figure}

In Figure \ref{fig:codifferentialGrid} we test numerically our discretization and compare it to the codifferential of \citet{Alexa2011}. We observe that the $L^2$ errors become constant. Concretely, for the jittered meshes with $r=0.4$, we get circa $5.69\cdot 10^{-2}$ for $\delta\beta$ and $2.82\cdot 10^{-1}$ for $\beta\mathsf{AW}$, i.e., our approximation error is roughly $5\times$ smaller. We have seen this difference on more or less irregular planar and non--planar meshes for trigonometric, linear, and quadratic forms.

\subsection{The Laplace Operator}\label{subsec:DiscreteLaplacian}
The Laplace--deRham operator $\Delta$ takes differential $k$--forms to $k$--forms and is defined as $ \Delta=d\delta + \delta d$, where $\delta$ is the codifferential and $d$ is the exterior derivative. On 0--forms (functions), it simplifies to 
\[
 \Delta= \delta d.
\]
We define our discrete Laplacian in the same manner, using our codifferential. 

\begin{figure}[htb]
\centering
\includegraphics[width=0.24\linewidth]{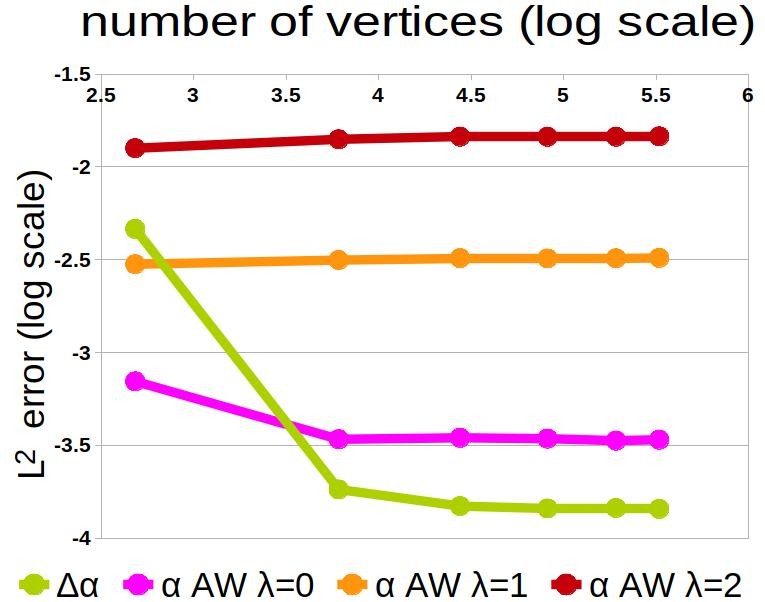}
\includegraphics[width=0.24\linewidth]{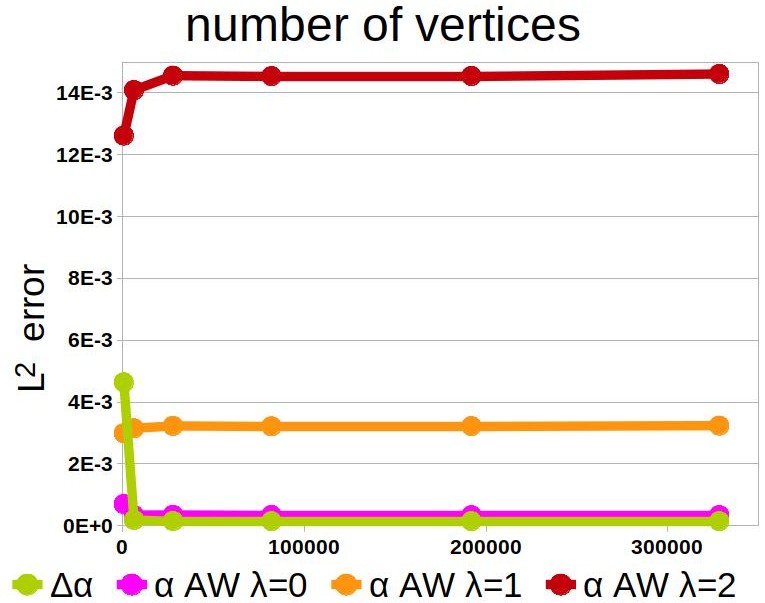}\hfill
\includegraphics[width=0.24\linewidth]{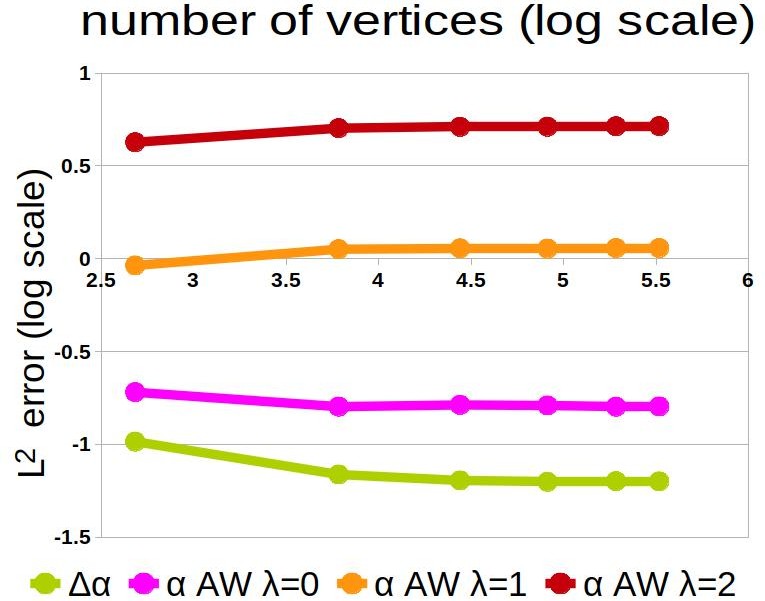}
\includegraphics[width=0.24\linewidth]{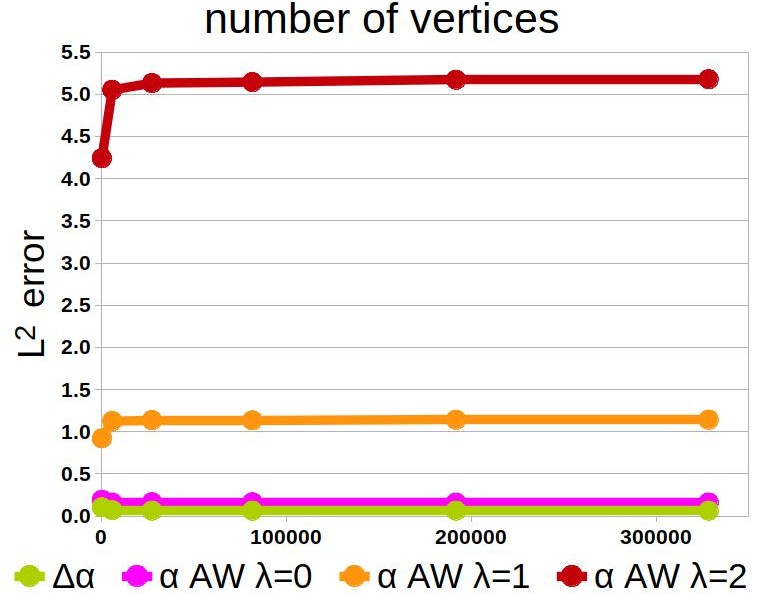}
\caption{The discrete Laplacian of a trigonometric 0--form $\alpha^0 = sin(x-1) - cos(2y)$ on two sets of planar quadrilateral jittered meshes with $r=0.01$ (far and center left) and $r=0.2$ (center and far right). We use the same set of meshes as in Figure \ref{fig:codifferentialGrid}. 
$\Delta\alpha$ denotes our Laplacian, $\alpha\,\mathsf{AW}\,\lambda=0$ the purely geometric, $\alpha\,\mathsf{AW}\,\lambda=1$ and $\alpha\,\mathsf{AW}\,\lambda=2$ the combinatorially enhanced Laplacians of \citet{Alexa2011}. The graphs in the center left and on the far right are in arithmetic scale. 
}\label{fig:LaplaceTrigGrid}
\end{figure}

Our Laplace operator is linearly precise, i.e., it is zero on linear 0--forms in plane. In Figure  \ref{fig:LaplaceTrigGrid} we depict the numerical behavior of our Laplacian on a trigonometric 0--form and compare it to the combinatorially enhanced (so called $\lambda$--simple choices with $\lambda=1$, $\lambda=2$) and purely geometric ($\lambda=0$) Laplacians of \citet{Alexa2011}. We note that all the $L^2$ errors become constant. Our and the purely geometric Laplacian give a better approximation to the analytical Laplacian than the combinatorially enhanced Laplacians. We have observed this pattern also on different quadratic and trigonometric 0--forms on more or less irregular polygonal meshes.

\section{Applications}\label{sec:Applications}
In this section we show some basic applications of our operators on general polygonal meshes.

\subsection{Implicit Mean Curvature Flow}\label{subsec:CurvatureFlow}
One of the widely used methods for smoothing a surface is the implicit mean curvature flow. If $f$ is a discrete 0--form representing vertex positions, then $\Delta f$ give us the direction and magnitude in which we should move each point in order to smooth the given mesh, see, e.g., \citep{CurvatureFlow}.

Let $f_0$ denote the initial state and $f_t$ the configuration after a mean curvature flow of some duration $t>0$. We employ the \textit{backward Euler scheme} to calculate $f_t$ by solving the linear system:
\[
 (I - t\Delta)f_t =f_0,
\]
where $I$ is the identity matrix. To solve this system, we use the \texttt{mldivide} algorithm of MATLAB.

In Figure \ref{fig:LaplaceFlowKitty} we show smoothing of general polygonal meshes and compare our method to the one of \citet{Alexa2011} with purely geometric Laplacians ($\lambda=0$) and combinatorially enhanced Laplacians ($\lambda=1$). After testing also other meshes and several other parameters $\lambda$, time steps, and number of iterations, we conclude that our results are visually comparable to theirs if $\lambda \in [1,2]$, and that our scheme does not create as many undesirable artifacts as theirs for $\lambda = 0$.

\subsection{Helmholtz--Hodge Decomposition}\label{subsec:HHD}
By the Hodge Decomposition Theorem \citep[Theorem 7.5.3]{Abraham}, if $M$ is a compact oriented Riemannian manifold without boundary and $\omega^k\in\Omega^k(M)$, then there exist uniquely determined forms $\alpha^{k-1}$, $\beta^{k+1}$, $\gamma^k$ ($\gamma$ harmonic, i.e., $\Delta\gamma = 0$) such that
\begin{equation}\label{eq:decomposition}
 \omega = d\alpha + \delta\beta + \gamma.
\end{equation}

If instead of forms, we think about a sufficiently smooth vector field $X = (\omega^1)^\sharp$, where $\sharp$ is the sharp operator, then an analogous Helmholtz theorem states that any vector field $X$ can be decomposed into an irrotational vector field (corresponding to $d\alpha$), a divergence--free component (analogous to $\delta\beta$), and a both irrotational and divergence--free vector field (corresponding to $\gamma$).
Thus the equation (\ref{eq:decomposition}) is also referred to as to Helmholtz--Hodge decomposition (HHD).

If $X$ is a divergence--free vector field (also known as solenoidal), we can find its two--component HHD, i.e., decompose $X$ into a rotational and irrotational part. In terms of differential forms, for $\omega^1 = X^\flat$ we get
\begin{equation}\label{eq:HHD2}
 \omega = \delta\beta + \gamma,
\end{equation}
where $\gamma$ is a harmonic 1--form and thus $\gamma^\sharp$ is an irrotational vector field, and $(\delta\beta)^\sharp$ is a rotational vector field. The two--component HHD is used for decomposition of vector fields of incompressible flows.

We use our codifferential operator to find our \textbf{discrete two--component Helmholtz--Hodge decomposition} as in equation (\ref{eq:HHD2}) by performing these steps:
\begin{enumerate}
 \item Discretize a given vector field $X$ with discrete flat operator (\ref{eq:flatOperator}) and define discrete 1--form $\omega^1 = X^\flat$.
 \item Find the 2--form $\beta$ by solving the equation $d\delta\beta = d\omega$.
 \item Set $\gamma = \omega - \delta\beta$.
\end{enumerate}

\begin{figure}[htb]
\centering\includegraphics[width=0.4\linewidth]{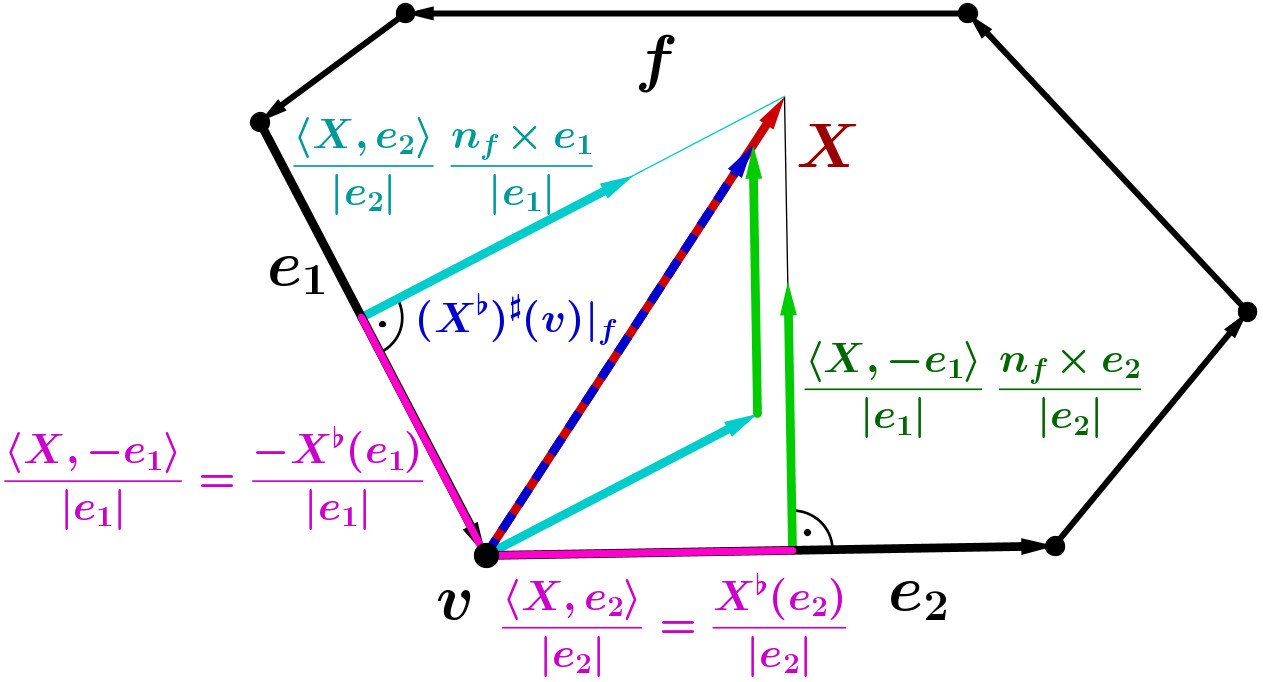}
\caption{The discrete sharp operator on a vertex $v$ restricted to a face $f$. Let $X^\flat$ be a 1--from computed by applying a discrete flat operator on a constant vector field $X$, then the orthogonal projection of $X$ on the unit direction vector of the edge $e_1$ equals $\frac{\langle X,e_1 \rangle}{|e_1|} = \frac{X^\flat(e_1)}{|e_1|}$. Similarly the orthogonal projection of $X$ on the unit direction vector of $e_2$ is $\frac{\langle X,e_2 \rangle}{|e_2|} = \frac{X^\flat(e_2)}{|e_2|}$.
Reconstructing the vector field $X$ from $X^\flat$, i.e., applying the sharp operator on $X^\flat$ as in equation (\ref{eq:sharpOperator}), yields vector
$(X^\flat)^\sharp|_f =  \frac{\langle X,e_2 \rangle}{|e_2|}\frac{n_f\times e_1}{|e_1|}
- \frac{\langle X,e_1 \rangle}{|e_1|}\frac{n_f\times e_2}{|e_2|}$
that has the same direction as $X$ and approximates its magnitude.
}\label{fig:SharpOperator}
\end{figure}

We can then map the discrete 1--forms $\delta\beta$ and $\gamma$ to discrete vector fields by applying \textbf{discrete sharp operator} $\sharp$ defined on an 1--form $\epsilon$ and per a vertex $v$ by:
\small
\begin{equation}\label{eq:sharpOperator}
\epsilon^\sharp (v) = \frac{1}{\rho(v)}\sum_{f \succ v}\bigg( \frac{\epsilon(e_2)}{|e_2|}\frac{n_f\times e_1}{|e_1|}
- \frac{\epsilon(e_1)}{|e_1|}\frac{n_f\times e_2}{|e_2|} \bigg),
\end{equation}
\normalsize
where $\rho(v)$ is the number of faces adjacent to $v$. Further $e_1,e_2\prec f$, $e_1$ is the edge which endpoint is $v$, $e_2$ is the edge with $v$ as the starting point, see Figure \ref{fig:SharpOperator}, and $n_f$ is a unit normal vector of the face $f=(v_0,\dots,v_{n-1})$ computed as:
\small
\[
 n_f = \frac{\hat{n}_f}{|\hat{n}_f|},\quad \hat{n}_f = \frac{1}{2}\sum_{j=0}^{n-1}(v_j\times v_{j+1}),\; \text{indices modulo } n.
\]
\normalsize

\begin{figure}[htb]\centering{
\includegraphics[width=0.24\linewidth]{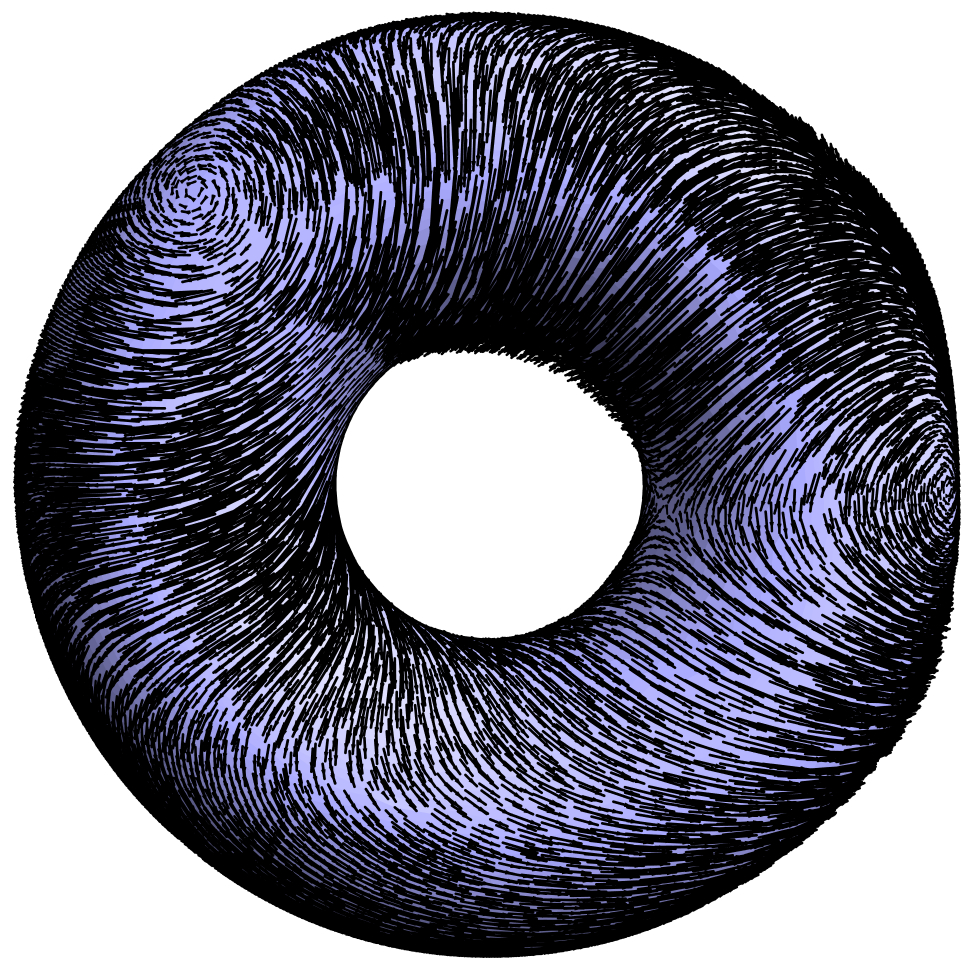}
\includegraphics[width=0.24\linewidth]{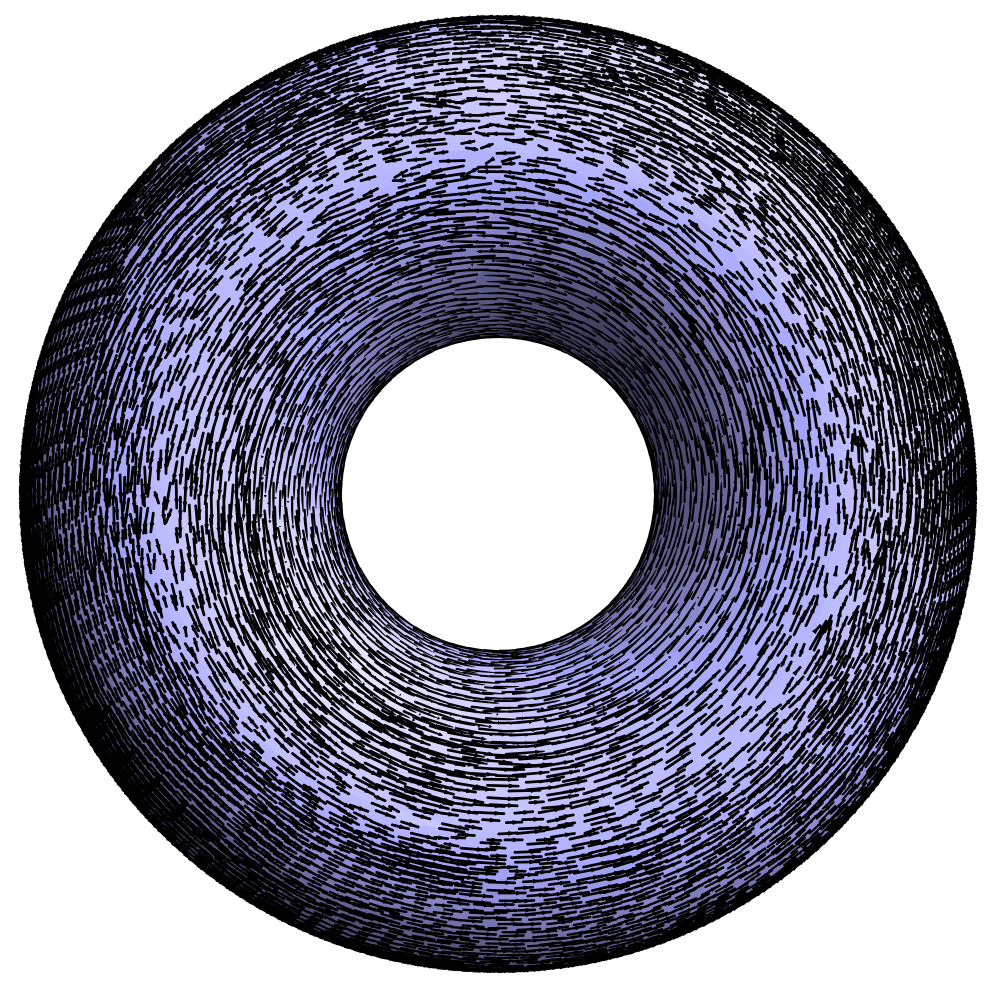}
\includegraphics[width=0.24\linewidth]{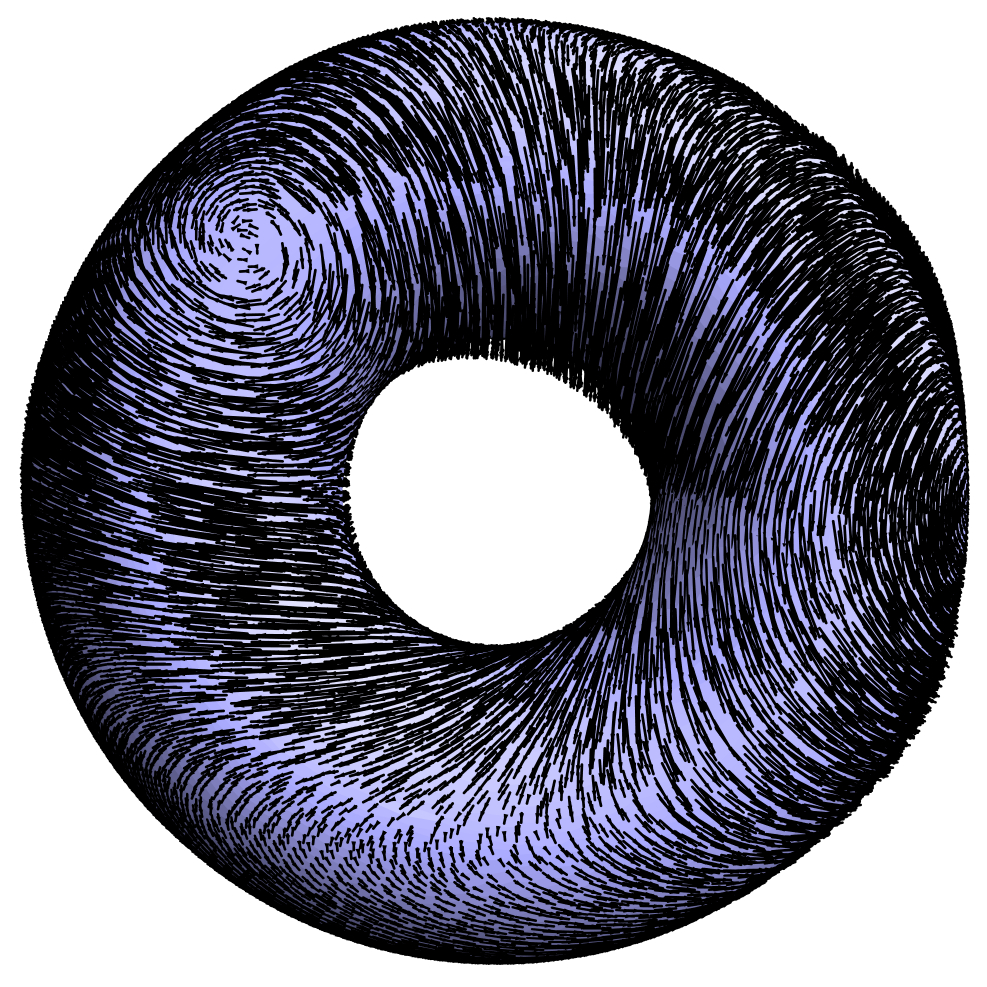}
\includegraphics[width=0.24\linewidth]{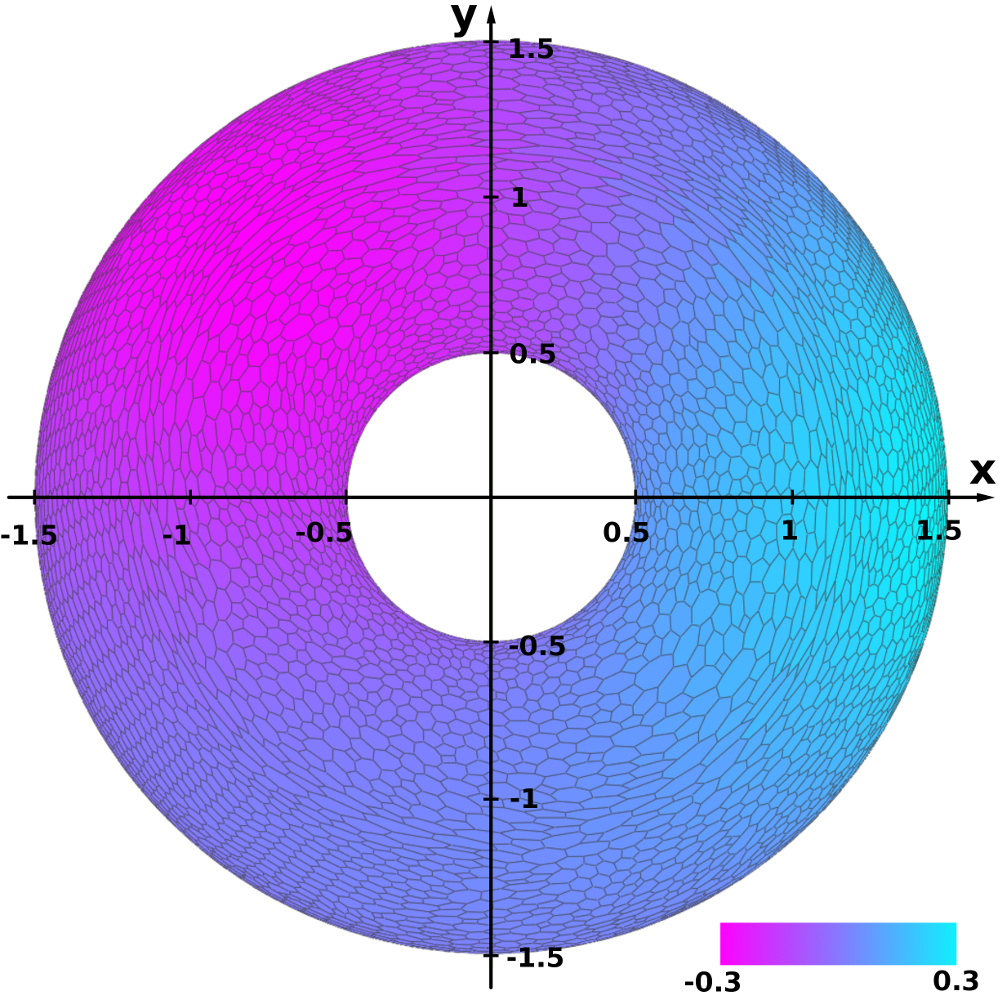}}
\caption{HHD of a solenoidal vector field $X$ on a torus centered at the origin (a general polygonal mesh with 20k vertices). The decomposed vector field is $X = X_H + X_R$, where
$X_H = (-y,x,0)$ is a harmonic field on the torus, and $X_R$ is a rotational vector field given by
$
X_R = \nabla(\exp(-(x-x_1)^2 - (y-y_1)^2 - (z-z_1)^2) - \exp(-(x-x_2)^2 - (y-y_2)^2 - (z-z_2)^2))\times n,
$
where $n$ is the unit normal vector of the torus.
We have chosen the center of CCW rotation $(x_1,y_1,z_1) = (\frac{3}{2},0,0)$, where the vector potential $\beta^\sharp$ reaches its maximum, and the center of CW rotation at $(x_2,y_2,z_2) = (\frac{-\sqrt{2}}{2},\frac{\sqrt{2}}{2},\frac{1}{2})$, where the vector potential $\beta^\sharp$ has its minimum (is negative). The vector field $X$ is shown on the far left.
Our discrete decomposition gives approximate expected results: the harmonic part $\gamma^\sharp$ calculated by our method is in the center left and the rotational part $(\delta\beta)^\sharp$ in the the center right. On the far right we visualize in pseudocolors the vector potential $\beta^\sharp$ computed by our algorithm.
}\label{fig:HHDTorus20k}
\end{figure}

\begin{figure}[htb]\centering
\includegraphics[width=0.29\linewidth]{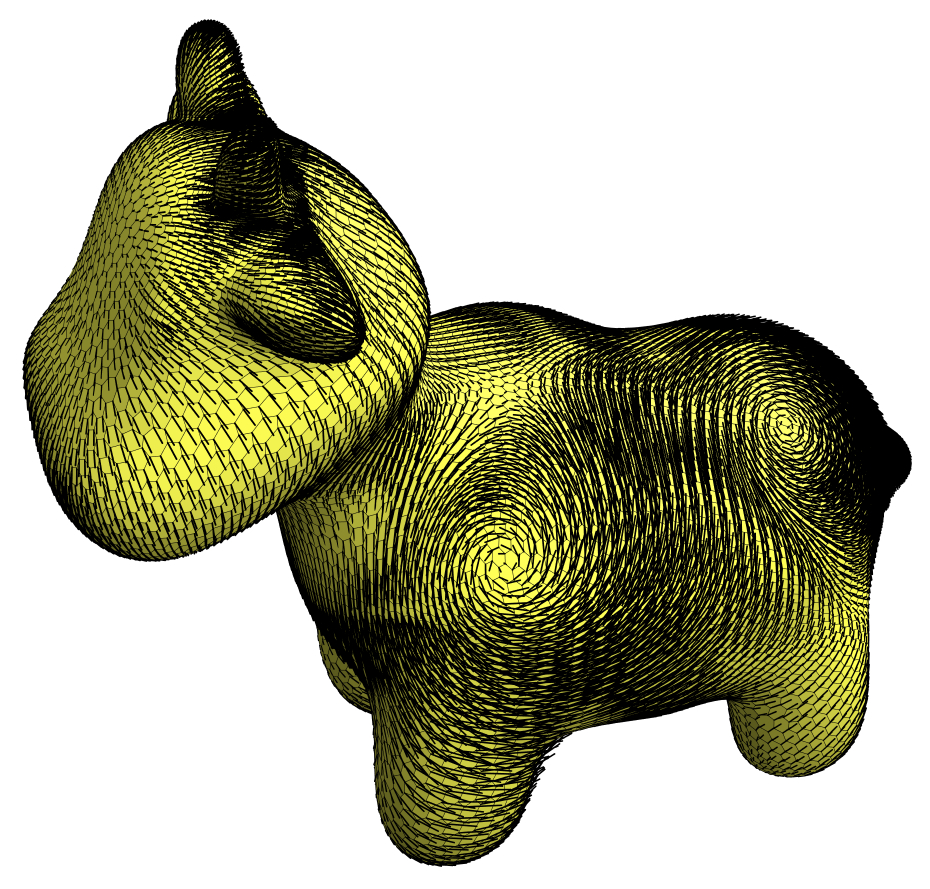}
\includegraphics[width=0.29\linewidth]{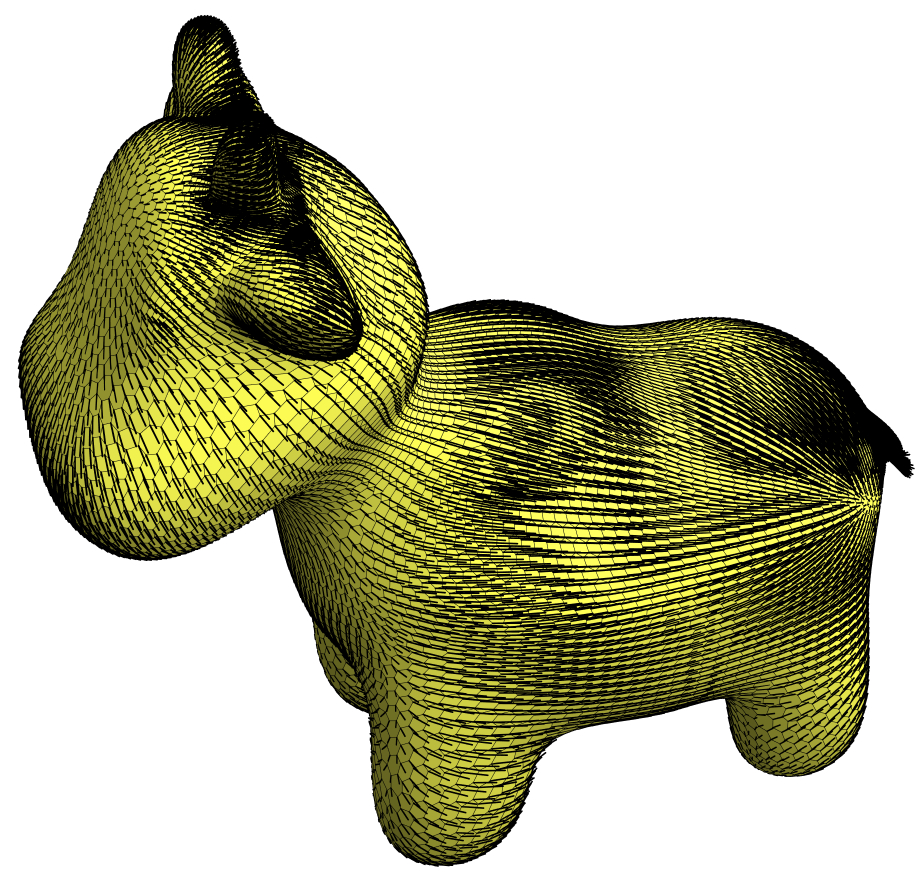}
\caption{HHD applied to remove the vortices of a vector field $\omega$ on a mesh of Spot (model created by Keenan Crane). On the left is the original vector field $\omega = \delta \beta + \phi$, on the right is its curl--free part $\phi$.}\label{fig:HHDCow}
\end{figure}

In Figure \ref{fig:HHDTorus20k} we give an example of our HHD of an incompressible vector field on a general polygonal mesh of a torus. In Figure \ref{fig:HHDCow} we then employ the HHD to remove vortices of an arbitrary vector field.

\subsection{Lie Advection}\label{subsec:LieAdvection}
The Lie derivative finds its application in dynamical systems. In computer graphics the Lie advection of differential forms (including scalar and vector fields) is used for tasks ranging from fluid flow simulation \citep{HOLA} to authalic parametrization of surfaces \citep{Zou}.

\begin{figure}[htb]
\begin{center}
\includegraphics[width=0.24\linewidth]{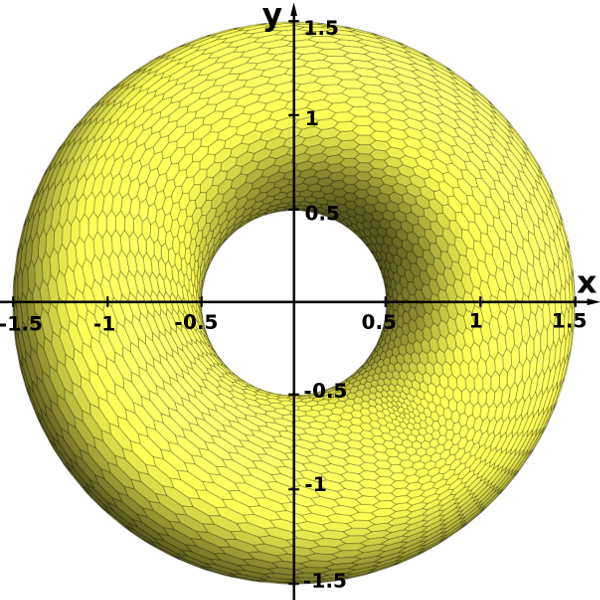}
\includegraphics[width=0.24\linewidth]{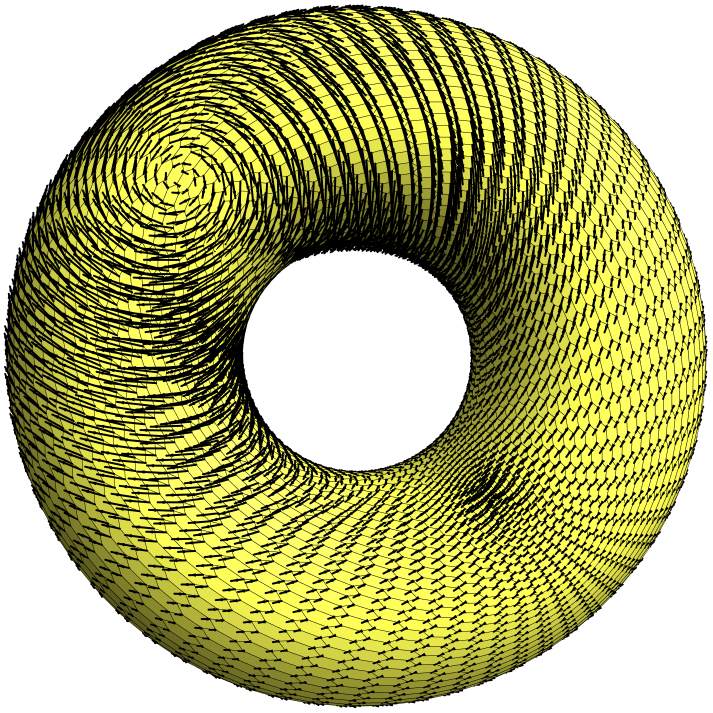} 
\includegraphics[width=0.24\linewidth]{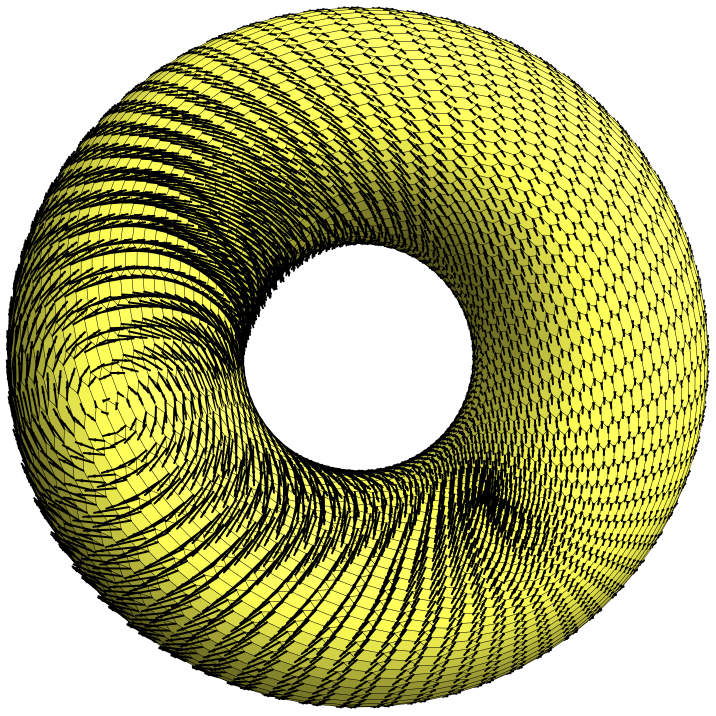} 
\includegraphics[width=0.24\linewidth]{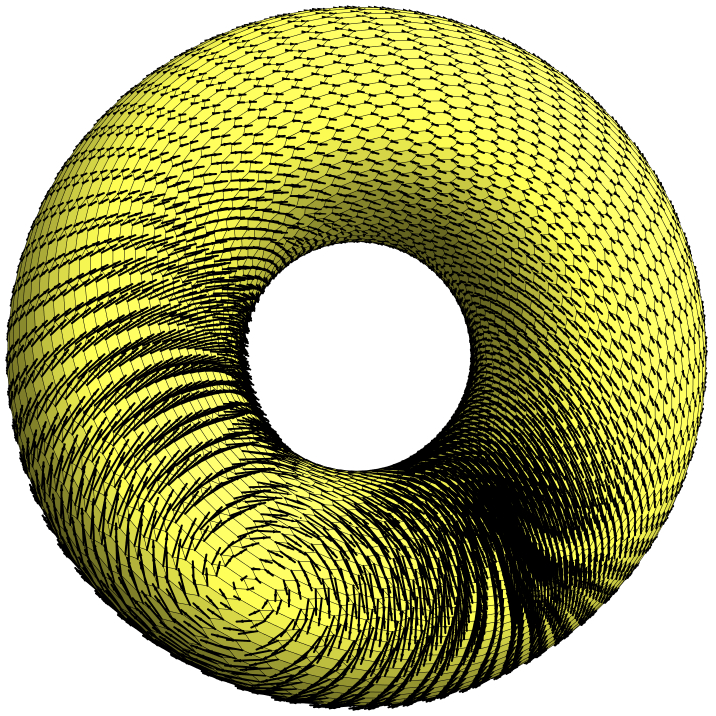} \\
\includegraphics[width=0.24\linewidth]{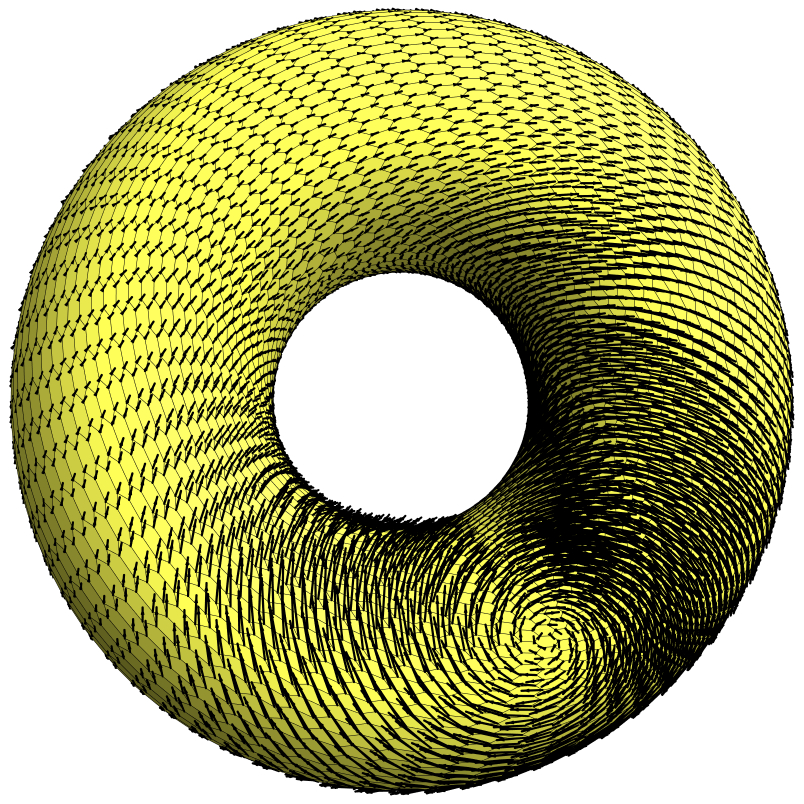}
\includegraphics[width=0.24\linewidth]{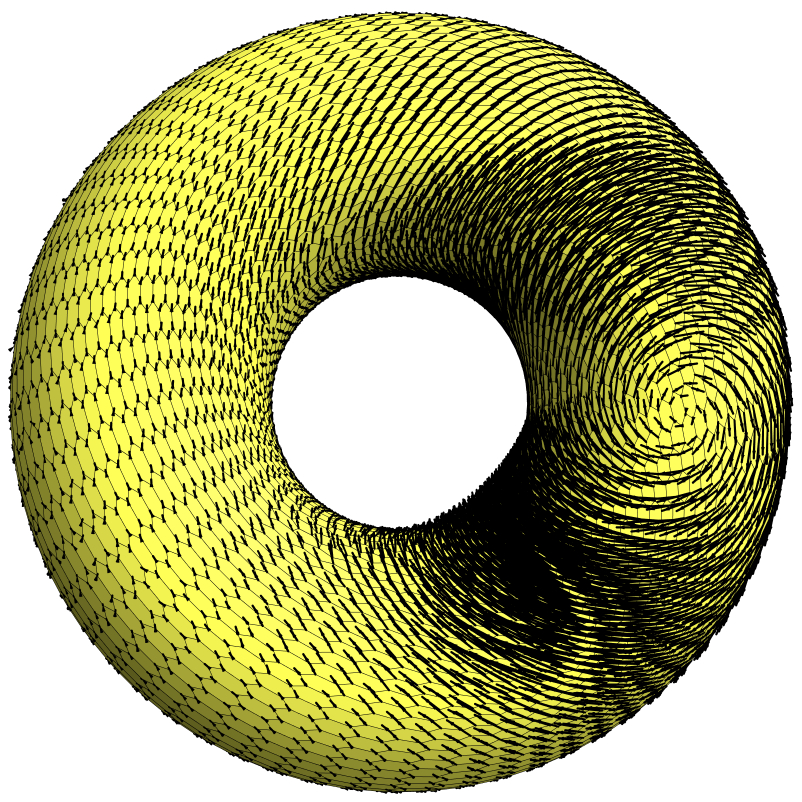} 
\includegraphics[width=0.24\linewidth]{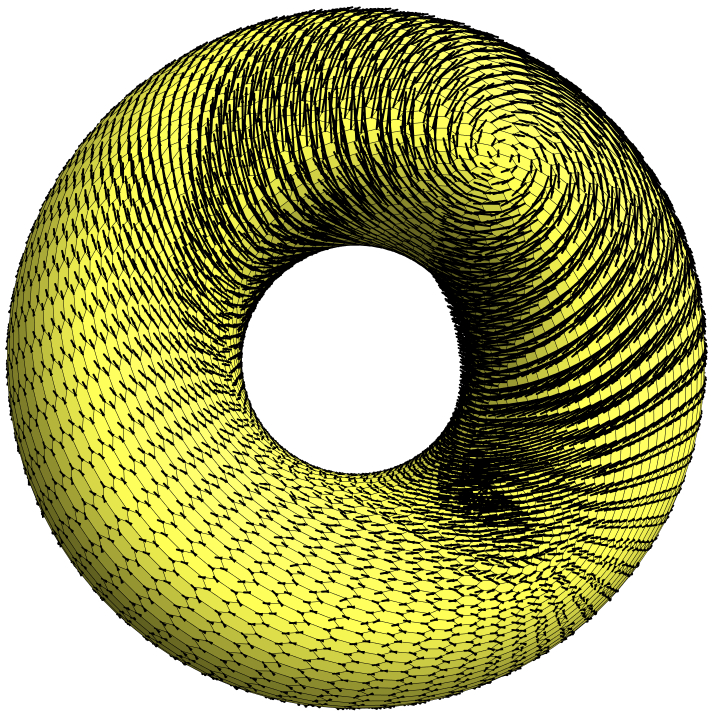} 
\includegraphics[width=0.24\linewidth]{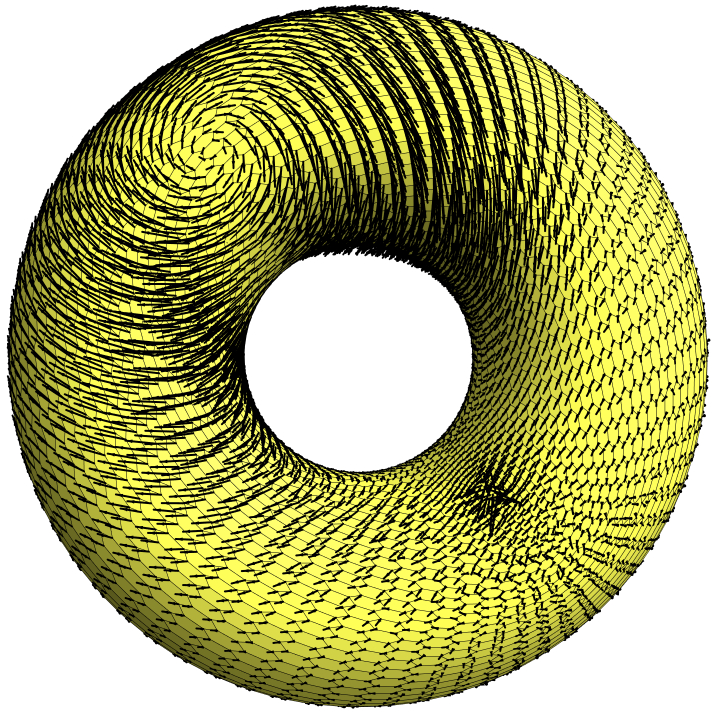} 
\end{center}
\caption{The Lie advection on a polygonal mesh on a torus (8k vertices), the mesh is shown in the top left corner. We advect an 1--form $\beta = Y^\flat$ along the flow of a tangent vector field $X=(-y,x,0)$, where $Y$ is given in equation (\ref{eq:YtoAdvect}). The second picture from the left depicts $\beta_0^\sharp$. We apply time steps of length $10^{-3}$. From left to right and top to bottom, we plot $Y =\beta^\sharp_k$ after 1000, 2000,\dots, 5000 iterations. At the bottom right we see $\beta_k^\sharp$ for $k=6283$. Because the domain is periodic, $\beta$ should be advected back to its original state $\beta_0$ after $2\pi\cdot 10^3 \approx 6283$ iterations. We can see that $\beta^\sharp_{6283}$ gets close to $\beta^\sharp_0$, but some undesirable artifacts appear -- especially at the small irregularity around the spot with coordinate position $(\frac{\sqrt{2}}{2},-\frac{\sqrt{2}}{2},\frac{1}{2})$.
}\label{fig:LieAdvectionTorus}
\end{figure}

In Figure \ref{fig:LieAdvectionTorus} we employ our Lie derivative to perform a simple discrete Lie advection of a tangent vector field $Y$ along a tangent vector field $X=(-y,x,0)$ on a torus azimuthally symmetric about the $z$ axis. $Y$ is a vorticial vector field given as
\begin{equation} \label{eq:YtoAdvect}
  Y = -\nabla\exp\bigg(-\Big(x+\frac{\sqrt{2}}{2}\Big)^2 - \Big(y-\frac{\sqrt{2}}{2}\Big)^2 - \Big(z-\frac{1}{2}\Big)^2 \bigg)\times n,
\end{equation}
where $n$ is the unit normal vector of the torus. To advect $Y$, we discretize it as a 1--form $\beta = Y^\flat$ and denote this initial state as $\beta_0$. We then iterate over our discrete solutions using a simple forward Euler method:
\begin{equation}\label{eq:LieAdvection}
\beta_{k+1} = \beta_k - t\Le_X\beta_k,\;k= 0,\dots,
\end{equation}
where $t$ is the time step, $k$ is the number of iterations, and each $\Le_X\beta_k$ is computed using our discrete Lie derivative (\ref{eq:discreteLie}). Note that the vector field $X$ is also discretized as a discrete 1--form.

The Lie derivative can be employed also for advection of a function by a vector field. In Figure \ref{fig:LieVaso} we advect a color function on a mesh of a vase.

\section{Discussion}\label{sec:Discussion}
Geometry processing with general polygonal meshes is a new developing area. We propose various discrete operators that act directly on meshes made of arbitrary polygons, possibly non--planar and non--convex, and thus open the possibility to perform many geometry processing tasks directly on these meshes.

Tangent vector fields on surfaces are used in many applications in computer graphics and other areas. We propose to represent vector fields as 1--forms and we provide methods for their design such as the Helmholtz--Hodge decomposition or the Lie advection. However, further applications are now ready to be explored, e.g., finding pairs of vector fields with zero Lie derivative for surface parameterization or employing Le advection of volume forms for area preserving parametrization.

Furthermore, we present a novel discrete Laplace operator that is numerically comparable to the purely geometric Laplacian of \citet{Alexa2011}, but results in a better mesh smoothing. On the other hand, our Laplacian gives a better numerical approximation to the analytically computed solutions than their combinatorially enhanced Laplacians, yet performs as well as theirs in smoothing of tested general polygonal meshes.



\end{document}